\documentclass[11pt,a4paper]{article}%
\usepackage{amssymb,amsmath, amsfonts}
\usepackage{graphicx,graphics}
\usepackage{mathtools}
\usepackage{bbold}
\usepackage[english]{babel}
\usepackage{epsfig,url}
\usepackage{bbm,theorem}
\usepackage{a4wide}
\usepackage{color}
\usepackage{enumerate}
\usepackage{calrsfs}
\usepackage[bookmarks=true]{hyperref}
\usepackage{bookmark}
\usepackage{amsmath}
\usepackage{amsfonts}
\usepackage{amssymb}
\usepackage{verbatim}
\usepackage{graphicx}%
\providecommand{\U}[1]{\protect\rule{.1in}{.1in}}
\providecommand{\U}[1]{\protect\rule{.1in}{.1in}}
\DeclareMathAlphabet{\pazocal}{OMS}{zplm}{m}{n}
\newtheorem{theorem}{Theorem}[section]

\newtheorem{lemma}[theorem]{Lemma}

\newtheorem{proposition}[theorem]{Proposition}
\newtheorem{remark}[theorem]{Remark}

\numberwithin{equation}{section}
\numberwithin{theorem}{section}
\newcommand{\qed}{\hfill$\Box$}
\newenvironment{proof}{\begin{trivlist}\item[]{\em Proof:}\/}{\qed\end{trivlist}}
\newenvironment{proofof}[1][Proof]{\noindent \textit{#1.} }{\ \qed}

\newcommand{\R}{{\mathbb R}}
\newcommand{\C}{{\mathbb C\hspace{0.05 ex}}}

\newcommand{\beq}{\begin{equation}}
\newcommand{\eeq}{\end{equation}}
\newcommand{\beqs}{\begin{eqnarray}}
\newcommand{\eeqs}{\end{eqnarray}}

\newcounter{jlisti}

\begin{document}

\title{Non-power law constant flux solutions for the Smoluchowski
coagulation equation}
\author{Marina A. Ferreira, Jani Lukkarinen, Alessia Nota, Juan J. L. Vel\'azquez}


\maketitle

\begin{abstract}
It is well known that for a large class of coagulation kernels, 
Smoluchowski coagulation equations have particular power law solutions  which yield a constant flux of mass along all  scales
of the system. In this paper, we prove that for some choices of the coagulation
kernels there are solutions   with a constant flux of mass along all  scales which
are not power laws. The result is proved by means of a bifurcation argument.

\end{abstract}


\tableofcontents


\section{Introduction}

\subsection{Motivation and general background of the problem}

In this paper, we study a particular type of stationary solutions of the
classical Smoluchowski coagulation equation
\begin{equation}
\partial_{t}f\left(  x,t\right)  =\mathbb{K}\left[  f\right]  \left(
x,t\right)
\end{equation}
where
\begin{align}
\mathbb{K}\left[  f\right]  \left(  x,t\right)  :=\frac{1}{2}\int_{0}%
^{x}K\left(  x-y,y\right)  f\left(  x-y,t\right)  f\left(  y,t\right)
dy-\nonumber\\
\int_{0}^{\infty}K\left(  x,y\right)  f\left(  x,t\right)  f\left(
y,t\right)  dy\ \ ,\ \ x>0\ \ ,\ \ t\geq0. \label{A1E1}%
\end{align}

The function $f$ describes the concentration of particles in the space of
particle volumes. The collision kernel $K$ is non-negative and symmetric, i.e.,
$K\left(  x,y\right)  =K\left(  y,x\right)  $ for each $x,y\in\mathbb{R}_{+}.$
This kernel contains information about the specific mechanism responsible for the
coagulation of the particles. See \cite{Fried, V} where coagulation kernels
have been obtained in the context of atmospheric aerosols. For a rigorous connection between coagulation kernels and dynamics of interacting particle
systems, see e.g. \cite{GKO, HR, LN, NoV, YRH}, and for dynamics of
stochastic processes in graphs, see e.g. \cite{A99, Du}.
Several of the collision kernels arising in the applications of the
Smoluchowski coagulation equation are homogeneous, i.e., there exists
$\gamma\in\mathbb{R}$ such that
\begin{equation}
K\left(  \lambda x,\lambda y\right)  =\lambda^{\gamma}K\left(  x,y\right)
\ \ \text{for any }x,y\in\mathbb{R}_{+}\text{ and }\lambda>0\, .\label{A1E2}%
\end{equation}

Equation (\ref{A1E1}) can be written in a more convenient form in order to
describe the transfer of  volume of the clusters from smaller to larger
values.  To this end, we set $g\left(  x,t\right)  =xf\left(  x,t\right)  .$ Then
(\ref{A1E1}) can be formally rewritten as
\begin{equation}
\partial_{t}\left(  g\left(  x,t\right)  \right)  +\partial_{x}\left(
J\left(  x;f\left(  \cdot,t\right)  \right)  \right)  =0 \label{S1E0}%
\end{equation}
where
\begin{equation}
J\left(  x;f\right)  =\int_{0}^{x}dy\int_{x-y}^{\infty}dzK\left(  y,z\right)
yf\left(  y\right)  f\left(  z\right)  . \label{S1E1}%
\end{equation}
Equation (\ref{S1E0}) shows that the coagulation process described by
Smoluchowski equation can be reinterpreted in terms of the flux of volume in
the space of cluster volumes. The total flux of volume passing from the region
of volumes smaller than $x$ to the region of volumes larger than $x$ is given
by (\ref{S1E1}).

A particular class of solutions of Smoluchowski coagulation equation are the
particle distributions $f$ for which the flux of particles at each particular
value of $x$ is a constant (and therefore independent of $t$ and $x$). Then,
\begin{equation}
J\left(  x;f\right)  =J_{0}>0\text{ \ for each }x>0\,. \label{A1E3}%
\end{equation}
In this paper, we will consider solutions to \eqref{A1E3} for the above mentioned
class of homogeneous kernels. Due to the homogeneity of the kernel, $K$ can then be written in
the form
\begin{equation}
K\left(  x,y\right)  =\left(  x+y\right)  ^{\gamma}\Phi\left(  \frac{x}%
{x+y}\right)  \label{S1E3}%
\end{equation}
where%
\begin{equation}
K\left(  s,1-s\right)  =\Phi\left(  s\right)  =\Phi\left(  1-s\right)
\ \ ,\ \ 0<s<1\,. \label{S1E4}%
\end{equation}
Notice that the last identity follows from the symmetry of the kernel $K$.

Besides the homogeneity of the kernel, the main property characterizing the
coagulation mechanism associated to a given coagulation kernel is the
asymptotic behaviour of the function $\Phi\left(  s\right)  $ as
$s\rightarrow0^{+}$ (or, equivalently, for $s\rightarrow1^{-}$). A large class of
kernels relevant to the applications of the Smoluchowski equation allows finding 
$0<c_{1}\leq c_{2}<\infty$ and $p\in\mathbb{R}$ such that 
\begin{align}
c_{1}s^{-p}\left(  1-s\right)  ^{-p}\leq\Phi\left(  s\right)  \leq c_{2}%
s^{-p}\left(  1-s\right)  ^{-p}\ \ \text{for }s\in\left(  0,1\right)\,.
\ \label{A1E6}%
\end{align}
In addition, we will show that it is possible to assume without loss of generality  that
\begin{equation}
\label{eq:cond_0}\gamma+2p\geq0 .
\end{equation}
We will assume here the above positivity, and the following upper bound
\begin{equation}
\gamma+2p<1 \ . \label{S1E7}%
\end{equation}

In order to justify that \eqref{eq:cond_0} can be assumed without loss of generality, we just notice that if $\gamma+2p<0$ then \eqref{S1E3} and \eqref{A1E6} hold with $p$ replaced by $\tilde p := -(\gamma+p)$. Using the fact that $\gamma+2\tilde p = -(\gamma+2p) \geq 0$ it follows that \eqref{A1E6} and \eqref{eq:cond_0} hold with $p$ replaced by $\tilde p$.

The class of kernels for which the conditions \eqref{S1E3}-\eqref{S1E7} hold,
covers all the kernels satisfying
\[
c_{1} (x^{\alpha}y^{\beta}+ x^{\beta}y^{\alpha}) \leq K(x,y) \leq c_{2}
(x^{\alpha}y^{\beta}+ x^{\beta}y^{\alpha}) \quad\text{ with } 0\leq
\alpha-\beta< 1\,,
\]
as it may be seen by choosing $\gamma=\alpha+\beta$ and $p = -\beta$.
Note that then, indeed, $0\le \gamma+2 p <1$.

Due to the homogeneity of the kernel $K$ we can expect the functional
$J\left(  x;f\right)  $ to be constant if%
\begin{equation}
f\left(  x\right)  =\frac{C_{0}}{x^{\frac{3+\gamma}{2}}},\ C_{0}>0,
\label{A1E7}%
\end{equation}
assuming that the integral in (\ref{S1E1}) is defined. It turns out that the
integral in (\ref{S1E1}) is finite with $f$ as in (\ref{A1E7}) if and only if
\eqref{S1E7} holds.
Indeed, using \eqref{A1E7} in \eqref{S1E1} and the upper bound in \eqref{A1E6},
we have
\begin{align*}
J\left(  x;f\right)   &  = \int_{0}^{1}dy\int_{1-y}^{\infty}dz\left(
y+z\right)  ^{\gamma}\Phi\left(  \frac{y}{y+z}\right)  \frac{y}{y^{\frac
{3+\gamma}{2}}}\frac{1}{z^{\frac{3+\gamma}{2}}}\\
&  \leq c_{2}\int_{0}^{1}dy\int_{1-y}^{\infty}dz\left(  y+z\right)  ^{\gamma
}\left(  \frac{y}{y+z}\right)  ^{-p}\left(  \frac{z}{y+z}\right)  ^{-p}%
\frac{y}{y^{\frac{3+\gamma}{2}}}\frac{1}{z^{\frac{3+\gamma}{2}}}.
\end{align*}
We now decompose the integral above in the contributions of the regions where
$y<z$ and 
$y\geq z$. We can then estimate $J\left(  x;f\right)  $ by
\begin{align}
J\left(  x;f\right)  \leq &  \ C\int_{0}^{1}dy\int_{\frac{1}{2}}^{\infty
}dz\left(  y+z\right)  ^{\gamma}\left(  \frac{y}{y+z}\right)  ^{-p}\left(
\frac{z}{y+z}\right)  ^{-p}\frac{y}{y^{\frac{3+\gamma}{2}}}\frac{1}%
{z^{\frac{3+\gamma}{2}}}+\nonumber\\
&  +C\int_{\frac{1}{2}}^{1}dy\int_{1-y}^{y}dz\left(  y+z\right)  ^{\gamma
}\left(  \frac{y}{y+z}\right)  ^{-p}\left(  \frac{z}{y+z}\right)  ^{-p}%
\frac{y}{y^{\frac{3+\gamma}{2}}}\frac{1}{z^{\frac{3+\gamma}{2}}}=: J_{1}+J_{2}
.
\end{align}
A straightforward computation leads to
\begin{align*}
&  J_{1}\leq C\int_{0}^{1}\frac{dy}{y^{\frac{1+\gamma}{2}+p}}\int_{\frac{1}%
{2}}^{\infty}dz\frac{\left(  z\right)  ^{\frac{\gamma}{2}+p}}{z^{\frac{3}{2}}}
<\infty\quad\text{if } \ \gamma+2p<1\,,\\
&  J_{2}\leq C\int_{\frac{1}{2}}^{1}dy\int_{1-y}^{y}\frac{dz}{z^{\frac
{3+\gamma}{2}+p}}\leq C\int_{0}^{1}dy\int_{1-y}^{1}\frac{dz}{z^{\frac
{3+\gamma}{2}+p}}\leq C\int_{0}^{1}\frac{dz}{z^{\frac{1+\gamma}{2}+p}}%
<\infty\quad\text{if } \ \gamma+2p<1\,.
\end{align*}
Hence, $J\left(  x;f\right)  <\infty$ if $\gamma+2p<1$.
\medskip


It has been seen in \cite{FLNV3} that (\ref{S1E7}) is a necessary and
sufficient condition for the existence of solutions of (\ref{A1E3}) including  non-homogeneous kernels. A similar
result for a slightly more restrictive class of coagulation kernels has been
found in \cite{FLNV1}.

Notice that if (\ref{S1E7}) holds, a solution of the equation (\ref{A1E3}) is
given by (\ref{A1E7}) with $C_{0}$ given by
\begin{equation}
C_{0}=b_{\Phi,\gamma}\sqrt{J_{0}}\ \ \text{where}\ \ b_{\Phi,\gamma}=\left(
\int_{0}^{1}dy\int_{1-y}^{\infty}dz\frac{\left(  y+z\right)  ^{\gamma}%
}{y^{\frac{\gamma+1}{2}}z^{\frac{\gamma+3}{2}}}\Phi\left(  \frac{y}%
{y+z}\right)  \right)  ^{-1/2}\ . \label{S1E8}%
\end{equation}
Notice that $b_{\Phi,\gamma}<\infty$, due to (\ref{S1E7}), since the integrand
is bounded by $\frac{z^{\frac{\gamma}{2}+p-\frac{3}{2}}} {y^{\frac{\gamma
+1}{2}+p}} \chi_{\{z\geq\frac1 2\}} \chi_{\{0\leq y\leq1\}}$ for $y\leq z $
and by $\frac{1}{z^{\frac{\gamma}{2}+p+\frac{3}{2}}}
\chi_{\{\max\{1-z,\frac1 2\} \leq y \leq1\}} $ for $y> z $.

Our goal is to prove the existence of other solutions of (\ref{A1E3}) which
are not power laws. These non-power law solutions will be shown to exist for a
large class of homogeneous kernels with homogeneity $\gamma$.

The power law solutions (\ref{A1E7}) are well known and they have been
extensively used in the analysis of coagulation problems (cf.~\cite{TIN}, as
well as \cite{CRZ} where an argument analogous to the one used in the
framework of wave turbulence introduced in \cite{ZF} has been applied to
coagulation equations). On the other hand, it has been proved in \cite{FLNV1}
that the stationary solutions of the coagulation equation with injection, i.e.,
the solutions of the equation
\begin{equation}
\mathbb{K}\left[  f\right]  \left(  x,t\right)  +\eta\left(  x\right)  =0
\label{StatCoag}%
\end{equation}
where $\eta$ is compactly supported in $\left(  0,\infty\right)  $, behave as a constant flux solution  as
$x\rightarrow\infty$, i.e., as a solution of
(\ref{S1E1}), (\ref{A1E3}). We remark that in the case of the constant kernel
this result has been considered in \cite{Dub}, including the convergence to equilibrium.

It is natural to ask if the only possible solutions of (\ref{S1E1}),
(\ref{A1E3}) are the power laws (\ref{A1E7}). In this paper, we prove that this
is not the case. More precisely, we will prove that there exist kernels $K$
satisfying (\ref{S1E3})--(\ref{S1E7}) for which, in addition to
the power law solutions (\ref{A1E7}), there are other constant flux solutions,
i.e., solutions of (\ref{S1E1}), (\ref{A1E3}) which have a form different from
(\ref{A1E7}).

\begin{figure}[ptb]
\centering
\includegraphics[scale=0.6]{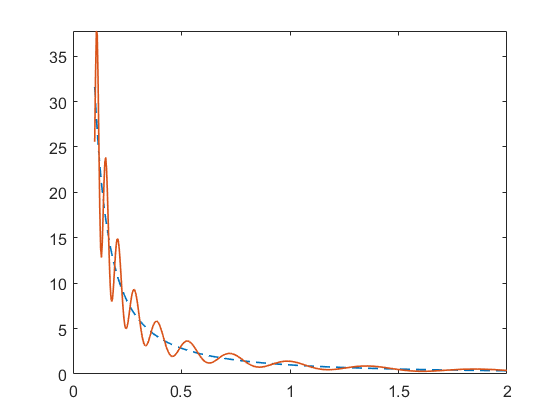}\caption{Representation of a power
law solution \eqref{A1E7} (dashed line) and a non-power law solution as in \eqref{eq:f_cos}.
}
\label{fig:1}
\end{figure}

The class of non-power law solutions of (\ref{S1E1}), (\ref{A1E3}) will be
obtained by means of a bifurcation argument. Specifically,
we will show in this paper (cf.\ Theorem \ref{main}) that there
are kernels 
$K$ 
satisfying  (\ref{S1E3})--(\ref{A1E6}) for which there exist
solutions of (\ref{S1E1}), (\ref{A1E3}) which have approximately the form%
\begin{equation}\label{eq:f_cos}
f\left(  x\right)  \simeq \frac{C_{0}}{x^{\frac{\gamma+3}{2}}}\left[  1+\varepsilon
\cos\left(  \alpha\log\left(  x\right)  \right)  \right]
\end{equation}
for some suitable $\alpha>0$ and a small $\varepsilon\neq 0$. 
We have illustrated one such solutions in Figure  \ref{fig:1}.


Actually, we will prove that the bifurcation of non-power law solutions to
(\ref{S1E1}), (\ref{A1E3}) can be obtained for kernels $K$ satisfying not only
(\ref{A1E6}) but also the more restrictive condition $\Phi\left(  s\right)  \sim as^{-p}\text{ as }s\rightarrow0^{+}\ \ \text{with
}p\in\mathbb{R}$, with $a>0$.
Note that rescaling the kernel and the flux we can also take $a=1$, without loss of generality. Therefore, we have
\begin{equation}
\Phi\left(  s\right)  \sim s^{-p}\text{ as }s\rightarrow0^{+}\ \ \text{with
}p\in\mathbb{R}. \label{S1E5}%
\end{equation}
Then, due to (\ref{S1E4}), we have also
\begin{equation}
\Phi\left(  s\right)  \sim  \left(  1-s\right)  ^{-p}\text{ as }s\rightarrow
1^{-}
\,. \label{S1E6}%
\end{equation}

Solutions of the Smoluchowski coagulation equation, as well as more general
coagulation-fragmentation models, which exhibit oscillations in the volume
variable $x$ have been found for several choices of the coagulation kernel
$K.$ In \cite{HNV}, it has been found that for coagulation $K\left(  x,y\right)  $
kernels strongly concentrated along the line $x\simeq y$ the self-similar
solutions yielding the asymptotic behaviour of the solutions of the
coagulation equation (without particle injection) develop an oscillatory
behaviour. Asymptotic expansions of these oscillatory behaviours have been
obtained in \cite{McNV} and \cite{NV}. The papers \cite{BNV2, BNV1} provide a
rigorous construction of a class of solutions of the coagulation-fragmentation
equation with kernels $K\left(  x,y\right)  $ concentrated near the line $x=y$
for which the solutions converge asymptotically to a sequence of Dirac masses.
In the particular case of the so-called diagonal kernels $K\left(  x,y\right)
=x^{\gamma+1}\delta\left(  x-y\right)  $, a full description of the long time
asymptotics of the solutions of the Smoluchowski coagulation equation has been
obtained in \cite{LNV}. It turns out that the solutions to (\ref{A1E1}) with
the kernel $K\left(  x,y\right)  =x^{\gamma+1}\delta\left(  x-y\right)  $ and
$\gamma<1,$ exhibit oscillations for most of the initial data both in time as well as 
in the volume variable $x.$

In all the examples described above, the results have been obtained for
coagulation or coagulation-fragmentation equations without injection of monomers. Ref.\
\cite{MSST} concerns a discrete coagulation equation for which injection of clusters
with two very different sizes takes place, namely monomers and clusters with
size $I_{s}\gg 1$. In the situation considered in \cite{MSST}, it 
would be natural to expect the solutions to behave for long times as
stationary solutions which follow a power law for  large cluster sizes. 
Indeed, this  turns
out to be the observed behaviour for large times. However, due to the
large difference of sizes between the clusters of size $I_{s}$ and the
monomers, the stationary solutions oscillate in the variable $x$ in a
large interval of cluster sizes $x$ until eventually the oscillations are damped
and the stationary solution finally approaches a power law for large values of $x.$
Contrarily, in the stationary solutions that we construct here, the
oscillations are present for all sizes and they result from properties of the coagulation
kernel rather than from the source term.

The present work does not concern oscillations in time that have been obtained
in the literature using both numerical and analytical methods (including
rigorous results), cf.\ \cite{ MKSTB,NPSV21, PV20}.

As remarked earlier, the condition (\ref{S1E7}) is a necessary and sufficient
condition for the existence of stationary solutions to (\ref{StatCoag}), i.e.,  with a source term. This problem has been
considered in the case of the discrete coagulation equation for kernels of the
form $K\left(  x,y\right)  =x^{\gamma+p}y^{-p}+y^{\gamma+p}x^{-p}$ and if the source term
$\eta\left(  x\right)  $ is a Kronecker delta at the monomers.  
In \cite{H},  formal asymptotic formulas for the concentration of large cluster sizes are derived using generating functions. These results indicate that solutions to  (\ref{StatCoag})  exist only if the condition (\ref{S1E7}) holds.
 More recently, in the case of both discrete and continuous coagulation equations
and general coagulation kernels, it has been rigorously proved in \cite{FLNV1}
that (\ref{S1E7}) is a necessary and sufficient condition for the existence of
solutions of (\ref{StatCoag}).  In fact, the condition (\ref{S1E7}) ensures the
so-called \textit{locality property} for the class of equations
(\ref{StatCoag}):
the most
relevant collisions are those between particles with comparable sizes and not
the collisions between particles with very different sizes  (cf.\ \cite{KC}). A similar property
has been extensively used in the study of the class of kinetic equations
arising in the theory of Wave Turbulence (cf.\ \cite{Zakh}).

\subsection{Main result of the paper}

In order to formulate the main result of this paper, it is convenient to
introduce some notation to characterize the class of admissible kernels. We
will denote as $\mathcal{K}_{\gamma,p}$ the class of continuous functions
$K\in C\big( \left(  0,\infty\right)  ^{2}\big) $ of the form (\ref{S1E3}),
(\ref{S1E4}) where $\Phi\in C\left(  0,1\right)  $ satisfies $\Phi\left(
s\right)  >0$ for $s\in\left(  0,1\right)  $ and the limit $a:=\lim
_{s\rightarrow0^{+}}\left[  s^{p}\Phi\left(  s\right)  \right]  $ exists and
is strictly positive. We endow $\mathcal{K}_{\gamma,p}$ with a structure of metric
space by means of the metric
\begin{equation}
\mathrm{dist} \left(  \Phi_{1},\Phi_{2}\right)  =\sup\left\{  s^{p}\left\vert
\Phi_{1}\left(  s\right)  -\Phi_{2}\left(  s\right)  \right\vert :s\in\left(
0,1\right)  \right\}  \label{A1E8}%
\end{equation}
We note that the metric space $\mathcal{K}_{\gamma,\lambda}$ is not complete
because the strict positivity of the kernels can be lost taking limits.

The following theorem presents the main result of this paper.

\begin{theorem}
\label{main}Let $J_{0}>0$. For each $\gamma,p\in\mathbb{R}$ satisfying \eqref{eq:cond_0} and
(\ref{S1E7}), there exists a one-parameter family
of kernels $K:I\rightarrow\mathcal{K}_{\gamma,p},$ with $I=\left(
-\delta,\delta\right)$, for some $\delta>0 $. The mapping $K$ is
continuous if $\mathcal{K}_{\gamma,p}$ is endowed with the topology generated
by the metric (\ref{A1E8}). Moreover, for each $\delta_{0} \in I \setminus
\left\{  0\right\}  $ there are at least two different solutions of
(\ref{S1E1}), (\ref{A1E3}). One of the solutions is given by (\ref{A1E7}) with
$C_{0}$ as in (\ref{S1E8}). The second solution $f_{\delta_{0} }$ has the
property that there exists $Q>1$ (independent of $\delta_{0} $) such that
$f_{\delta_{0} }\left(  Qx\right)  =Q^{-\frac{\gamma+3}{2}}f_{\delta_{0}
}\left(  x\right)  $ and the function $f_{\delta_{0} }$ is not a power law of
the form \eqref{A1E7} in the interval $\left(  1,Q\right)  .$
\end{theorem}

\begin{remark}
The kernels obtained in the proof of the Theorem are of the form (\ref{S1E3}), (\ref{S1E4}) with $\Phi\in
C^{\infty}\left(  0,1\right)$.  As an example with $\gamma=0$, we have plotted 
a representation of the non-power law solution constructed in the proof  in
Figure \ref{fig:1}.
\end{remark}

\subsection{Structure of the paper and main notation}

To quantify asymptotic properties of functions, we rely here on the following 
fairly standard notations.  We write ``$f\sim g$ as
$x\rightarrow\infty$'' to denote $\frac{f(x)}{g(x)}\rightarrow1$, as
$x\rightarrow\infty$. Moreover, given two functions $f,\ g$, we write ``$f\approx g$ in an interval $I\subset \R$'' if $\frac{g(x)}{2}\leq
f(x)\leq2g(x)$  for $x\in I\subset \R$, and the notation
``$f\ll g$'' is used  if the quotient $\frac{g(x)}{f(x)}$ can be made arbitrarily large for $x$ sufficiently large.   

The complex conjugate of $a\in \C$ is denoted by $\bar{a}$.
The indicator function of any set $A\subset{\mathbb{R}}$ will be denoted by
$\chi_{A}$.

The plan of the paper is the following. In Section \ref{ssec:2.1}, we
reformulate the problem (\ref{S1E1}), (\ref{A1E3}) using a more convenient set of
variables in such a way that power law solutions \eqref{A1E7} become a
constant solution. Then, in Section \ref{ssec:2.2}, we formulate the main properties of the linearized version of the problem around this power law  solution which are proved later in Section \ref{SectionLinear}.
 In Section
\ref{ssec:bifurcation}, we study the full nonlinear problem using a Hopf
bifurcation type of argument, concluding the proof of Theorem
\ref{main}. 

\section{Proof of the main result}


\label{sec:2}

In order to prove Theorem \ref{main}, it is convenient to reformulate the
problem (\ref{S1E1})--(\ref{A1E3}) using a different set of variables in which
the solutions (\ref{A1E7}) become constant. We will discuss later in this
section how to linearize the problem (\ref{S1E1}) around the power law
solution  or, equivalently, the reformulated problem around the constant
solution. The information obtained from the linearized problem will be used
later to prove a bifurcation result for the full nonlinear problem that will
imply Theorem \ref{main}.

\subsection{Reformulation of the problem}
\label{ssec:2.1} 

Here, we reformulate the problem (\ref{S1E1})--(\ref{A1E3}) so
that the solutions (\ref{A1E7}) become constant. Notice that (\ref{S1E1}%
)--(\ref{A1E3}) is invariant under a rescaling group. In the new set of
variables that we introduce in this section, the rescaling group becomes the
 group of translations. This will be convenient in order to bifurcate the non-constant flux solutions that we study in this paper. 
 
 We define a function
$H:\mathbb{R}\rightarrow\mathbb{R}$ such that
\begin{equation}
H\left(  X\right)  =x^{\frac{\gamma+3}{2}}f\left(  x\right)  \ \ ,\ \ x=e^{X}.
\label{S2E1}%
\end{equation}
Then, setting $y=e^{Y},\ z=e^{Z}$ we obtain
\[
K\left(  y,z\right)  yf\left(  y\right)  f\left(  z\right)  =\frac{1}%
{y^{\frac{1}{2}}z^{\frac{3}{2}}}W\left(  Y-Z\right)  H\left(  Y\right)
H\left(  Z\right)
\]
where
\[
W\left(  Y\right)  =\frac{1}{2}\left(  e^{\frac{Y}{2}}+e^{-\frac{Y}{2}%
}\right)  ^{\gamma}\left[  \Phi\left(  \frac{1}{1+e^{Y}}\right)  +\Phi\left(
\frac{1}{1+e^{-Y}}\right)  \right]  .
\]
Using (\ref{S1E4}), we find that $\Phi\left(  \frac{1}{1+e^{Y}}\right)
=\Phi\left(  \frac{1}{1+e^{-Y}}\right)  ,$ and thus
\begin{equation}
W\left(  Y\right)  =\left(  e^{\frac{Y}{2}}+e^{-\frac{Y}{2}}\right)  ^{\gamma
}\Phi\left(  \frac{1}{1+e^{Y}}\right)  =\left(  e^{\frac{Y}{2}}+e^{-\frac
{Y}{2}}\right)  ^{\gamma}\Phi\left(  \frac{1}{1+e^{-Y}}\right)  \ \ ,\ \ Y\in
\mathbb{R}\,. \label{S1E9}%
\end{equation}
Therefore, $W$ is symmetric, $W\left( - Y\right)  =W\left(  Y\right)  $.
Moreover, (\ref{S1E5}), (\ref{S1E6}) imply
\begin{equation}
W\left(  Y\right)  \sim  e^{\left(  \frac{\gamma}{2}+p\right)  \left\vert
Y\right\vert }\text{ as }\left\vert Y\right\vert \rightarrow\infty \,.\label{S2E3}%
\end{equation}

We further observe that, given any function $W$ satisfying $W\left(  Y\right)
=W\left(  -Y\right)  $ and (\ref{S2E3}), we can obtain a function $\Phi:\left(
0,1\right)  \rightarrow\mathbb{R}$ such that (\ref{S1E9}) holds. Indeed, using
the change of variable $\frac{1}{1+e^{Y}}=s$ in \eqref{S1E9}, we obtain
\begin{equation}
\Phi\left(  s\right)  =\frac{1}{\left(  \sqrt{\frac{1-s}{s}}+\sqrt{\frac
{s}{1-s}}\right)  ^{\gamma}}W\left(  \log\left(  \frac{1-s}{s}\right)
\right)  \ \ ,\ \ s\in\left(  0,1\right) \,. \label{S2E3a}%
\end{equation}
\begin{figure}[ptb]
\centering
\includegraphics[scale=0.6]{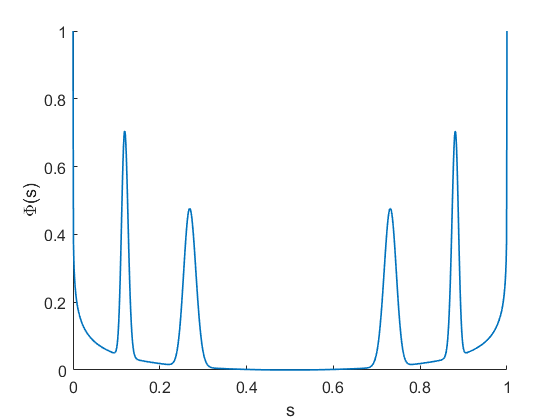}\caption{Representation of the
shape of a typical function $\Phi>0$ (cf.\ \eqref{S2E3a}) associated to the coagulation
kernel (cf.\ \eqref{S1E3}) considered in this paper, with $\gamma=0$ and
$p=0$. }%
\label{fig:phi}%
\end{figure}

A representation of a typical function $\Phi$ appearing later in the proof has been illustrated in Figure
\ref{fig:phi}. 

We can now reformulate the problem \eqref{S1E1}--\eqref{A1E3}
as
\begin{equation}
B\left(  H,H;W\right)  \left(  X\right)  =J_{0}\ \ ,\ \ X\in\mathbb{R}\,,
\label{S2E2}%
\end{equation}
where
\begin{equation}
B\left(  H_{1},H_{2};W\right)  \left(  X\right)  :=\int_{-\infty}^{X}%
dY\int_{X+\log\left(  1-e^{Y-X}\right)  }^{\infty}dZ\left[  e^{\frac{1}%
{2}\left(  Y-Z\right)  }W\left(  Y-Z\right)  \right]  H_{1}\left(  Y\right)
H_{2}\left(  Z\right) \ . \label{S2E4}%
\end{equation}
We will use repeatedly the following property of the operator $B.$
\begin{proposition}
For each $W\in L_{loc}^{\infty}\left(  \mathbb{R}\right)  $ satisfying
(\ref{S2E3}) the operator $B$ defined by means of (\ref{S2E4}) defines a
bilinear operator from $L^{\infty}\left(  \mathbb{R}\right)  \times L^{\infty
}\left(  \mathbb{R}\right)  $ to $\mathbb{R}$ which satisfies%
\begin{equation}
\left\Vert B\left(  H_{1},H_{2};W\right)  \left(  \cdot\right)  \right\Vert
_{L^{\infty}\left(  \mathbb{R}\right)  }\leq C_{W}\left\Vert H_{1}\left(
\cdot\right)  \right\Vert _{L^{\infty}\left(  \mathbb{R}\right)  }\left\Vert
H_{2}\left(  \cdot\right)  \right\Vert _{L^{\infty}\left(  \mathbb{R}\right)
} \label{S2E6}%
\end{equation}
where%
\begin{equation}
C_{W}=\int_{-\infty}^{0}dY\int_{\log\left(  1-e^{Y}\right)  }^{\infty
}dZ\left[  e^{\frac{1}{2}\left(  Y-Z\right)  }\left\vert W\left(  Y-Z\right)
\right\vert \right]  < \infty . \label{S2E6a}%
\end{equation}

\end{proposition}

\begin{proof}
From the definition of the bilinear operator $B$ (cf.~(\ref{S2E4})) we
immediately obtain the estimate
\begin{align}
\left\Vert B\left(  H_{1},H_{2};W\right)  \left(  \cdot\right)  \right\Vert
_{L^{\infty}\left(  \mathbb{R}\right)  }\leq &  \ \left\Vert H_{1}\left(
\cdot\right)  \right\Vert _{L^{\infty}\left(  \mathbb{R}\right)  }\left\Vert
H_{2}\left(  \cdot\right)  \right\Vert _{L^{\infty}\left(  \mathbb{R}\right)
}\nonumber\\
&  \ \sup_{X\in\mathbb{R}}\left(  \int_{-\infty}^{X}dY\int_{\log\left(
1-e^{Y-X}\right)  +X}^{\infty}dZ\left[  e^{\frac{1}{2}\left(  Y-Z\right)
}\left\vert W\left(  Y-Z\right)  \right\vert \right]  \right)  . \label{S2E5}%
\end{align}
Note that the integral on the right-hand side of this formula is independent
of $X$, hence
\begin{align*}
&  \sup_{X\in\mathbb{R}}\left(  \int_{-\infty}^{X}dY\int_{X+\log\left(
1-e^{Y-X}\right)  }^{\infty}dZ\left[  e^{\frac{1}{2}\left(  Y-Z\right)
}\left\vert W\left(  Y-Z\right)  \right\vert \right]  \right) \\
&  =\int_{-\infty}^{0}dY\int_{\log\left(  1-e^{Y}\right)  }^{\infty}dZ\left[
e^{\frac{1}{2}\left(  Y-Z\right)  }\left\vert W\left(  Y-Z\right)  \right\vert
\right]  \leq C\int_{-\infty}^{0}dY\int_{\log\left(  1-e^{Y}\right)  }%
^{\infty}dZ \left[  e^{\frac{1}{2}\left(  Y-Z\right)  }e^{\left(  \frac
{\gamma}{2}+p\right)  \left\vert Y-Z\right\vert }\right]
\end{align*}
where we used (\ref{S2E3}) in the last inequality. We now show that the
integral on the right hand side of the equation above is finite. Indeed,
\begin{align*}
&  \int_{-\infty}^{0}dY\int_{\log\left(  1-e^{Y}\right)  }^{\infty}dZ\left[
e^{\frac{1}{2}\left(  Y-Z\right)  }e^{\left(  \frac{\gamma}{2}+p\right)
\left\vert Y-Z\right\vert }\right] \\
&  =\int_{-\infty}^{0}dY\int_{\max\left\{  \log\left(  1-e^{Y}\right)
,Y\right\}  }^{\infty}dZ\left[  e^{\frac{1}{2}\left(  Y-Z\right)  }e^{\left(
\frac{\gamma}{2}+p\right)  \left(  Z-Y\right)  }\right] \\
&  \ \ \ +\int_{-\log(2)}^{0}dY\int_{\log\left(  1-e^{Y}\right)  }%
^{\max\left\{  \log\left(  1-e^{Y}\right)  ,Y\right\}  }dZ\left[  e^{\frac
{1}{2}\left(  Y-Z\right)  }e^{\left(  \frac{\gamma}{2}+p\right)  \left(
Y-Z\right)  }\right] \\
&  \leq\int_{-\infty}^{0}dY\int_{\max\left\{  Y,-\log\left(  2\right)
\right\}  }^{\infty}dZ\left[  e^{\left(  \frac{\gamma}{2}+p-\frac{1}%
{2}\right)  \left(  Z-Y\right)  }\right]  +\int_{-\log\left(  2\right)  }%
^{0}dY\int_{\log\left(  1-e^{Y}\right)  }^{Y}dZ\left[  e^{\left(  \frac{1}%
{2}+\frac{\gamma}{2}+p\right)  \left(  Y-Z\right)  }\right] \\
&  \leq\int_{-\infty}^{0}dY\int_{-\log\left(  2\right)  }^{\infty}dZ\left[
e^{\left(  \frac{\gamma}{2}+p-\frac{1}{2}\right)  \left(  Z-Y\right)
}\right]  +\int_{-\infty}^{0}dZ\int_{\max\left\{  \log\left(  1-e^{Z}\right)
,-\log\left(  2\right)  \right\}  }^{0}dY\left[  e^{\left(  \frac{1}{2}%
+\frac{\gamma}{2}+p\right)  \left(  Y-Z\right)  }\right] \\
&  \leq C_{1}\left(  \gamma,p\right)  +C\int_{-\infty}^{0}dZ\left[
e^{-\left(  \frac{1}{2}+\frac{\gamma}{2}+p\right)  Z}e^{Z}\right]  \leq
C_{1}\left(  \gamma,p\right)  +C_{2}\left(  \gamma,p\right)
\end{align*}
where $C_{1}\left(  \gamma,p\right)  ,\ C_{2}\left(  \gamma,p\right)  $ are
finite due to (\ref{S1E7}). Combining this with (\ref{S2E5}) we obtain
(\ref{S2E6}) with $C_{W}$ as in (\ref{S2E6a}).
\end{proof}

\medskip

In the variables \eqref{S2E1}, the constant flux solution
\eqref{A1E7},\ \eqref{S1E8} becomes the following solution of \eqref{S2E2}
\begin{equation}
H_{s}\left(  X\right)  =b_{\Phi,\gamma}\sqrt{J_{0}} \label{S2E7}%
\end{equation}
with $b_{\Phi,\gamma}$ as in (\ref{S1E8}). Notice for further reference that
\begin{equation}
B\left(  H_{s},H_{s};W\right)  \left(  X\right)  =J_{0} \ . \label{S2E7a}%
\end{equation}

Our goal is to prove the existence of solutions to (\ref{S2E2}) different from
(\ref{S2E7}) using a bifurcation argument. More precisely, we will obtain
solutions of (\ref{S2E2}) which are different from (\ref{S2E7}) but are close
to constant for some particular choices of kernel. To this end we first
consider a linearized version of (\ref{S2E2}).

\subsection{Linearized problem near the constant solutions.}

\label{ssec:2.2}

We first study the linearized problem obtained from (\ref{S2E2}) for small
perturbations of (\ref{S2E7}). To this end, we write
\begin{equation}
H\left(  X\right)  =H_{s}\left(  X\right)  \left[  1+\varphi\left(  X\right)
\right]  =b_{\Phi,\gamma}\sqrt{J_{0}}\left[  1+\varphi\left(  X\right)
\right]  . \label{S2E8}%
\end{equation}
Our goal is to obtain a class of kernels $W_{0}$ for which the linearized
problem has non trivial solutions. Plugging (\ref{S2E8}) into (\ref{S2E2}),
assuming that the kernel $W=W_{0}$, and neglecting quadratic terms in $\varphi$,
we obtain the linearized problem
\begin{equation}
\mathcal{L}\left(  \varphi;W_{0}\right)  \equiv B\left(  1,\varphi
;W_{0}\right)  +B\left(  \varphi,1;W_{0}\right)  =0 . \label{S2E9}%
\end{equation}
Using (\ref{S2E4}) we can rewrite (\ref{S2E9}) as
\begin{equation}
\mathcal{L}\left(  \varphi;W_{0}\right)  \left(  X\right)  =\int_{-\infty}%
^{X}dY\int_{X+\log\left(  1-e^{Y-X}\right)  }^{\infty}dZ\left[  e^{\frac{1}%
{2}\left(  Y-Z\right)  }W_{0}\left(  Y-Z\right)  \right]  \left(
\varphi\left(  Y\right)  +\varphi\left(  Z\right)  \right)  =0 . \label{S3E1}%
\end{equation}

Since we want to obtain solutions of the nonlinear problem (\ref{S2E2}) by
means of a perturbative argument, it is natural to look for bounded, non
trivial solutions of (\ref{S3E1}). Moreover, due to the fact that the operator
$\mathcal{L}\left(  \cdot;W_{0}\right)  $ commutes with the
group of translations, we look for solutions of (\ref{S3E1})
with the form $\varphi\left(  X\right)  =e^{ikX}$ for some $k\in\mathbb{R}$.
In order to have nonconstant solutions, we need $k\in\mathbb{R}\setminus
\left\{  0\right\}  .$  We  later prove the following result.
\begin{theorem}
\label{LinearzPb}For each $\gamma,\ p$ satisfying \eqref{eq:cond_0} and
\eqref{S1E7}, there exists a function $W_{0}=W_{0}\left(  X\right)  $,
$W_{0}\in H^{\infty}(D_{\beta})$ where $H^{\infty}(D_{\beta})$ denotes the
Hardy space (cf.\cite{Ru87}) of bounded, analytic functions in the open
domain
\begin{equation}
\label{def:Dbeta}D_{\beta}=\left\{  X\in\mathbb{C}:\left\vert \text{Im}\left(
X\right)  \right\vert < \left(  1+\beta\left\vert \text{Re}\left(  X\right)
\right\vert \right)  \right\}  \subset\mathbb{C}%
\end{equation}
for some $\beta>0$ (depending on $\gamma,p$). Moreover, $W_{0}$ satisfies

\begin{itemize}
\item[(i)] $W_{0}\left(  X\right)  =W_{0}\left(  -X\right)  $ for $X\in
D_{\beta}$

\item[(ii)] $W_{0}\left(  X\right)  >0$ for $X\in\mathbb{R}$

\item[(iii)] the following asymptotics (cf.\ (\ref{S2E3}))
\begin{equation}
W_{0}\left(  X\right)  \sim\exp\left(  \left(  \frac{\gamma}{2}+p\right)
\sqrt{X^{2}+1}\right)  \left[  1+O\left(  e^{-\kappa\left\vert X\right\vert
}\right)  \right]  \label{S3E2}%
\end{equation}
as $\left\vert X\right\vert =\sqrt{\left(  \text{Re}\left(  X\right)  \right)
^{2}+\left(  \text{Im}\left(  X\right)  \right)  ^{2}} \rightarrow\infty$,
$X\in D_{\beta}$ with $\kappa>0.$
\end{itemize}

Moreover, there exists a function $\Psi\left(  k;W_{0}\right)  $ analytic in
$D_{\beta}$ such that
\begin{equation}
\mathcal{L}\left(  e^{ik\cdot};W_{0}\right)  =\Psi\left(  k;W_{0}\right)
e^{ik\cdot}\ \ \ ,\ \ \ k\in\mathbb{R}. \label{S3E3a}%
\end{equation}
The function $\Psi\left(  k;W_{0}\right)  $ satisfies $\Psi\left(
-k;W_{0}\right)  =\overline{ \Psi\left(  k;W_{0}\right)  } $ and there exists
$k_{\ast}\in\mathbb{R},\ k_{\ast}>0$ such that
\begin{equation}
\Psi\left(  k_{\ast};W_{0}\right)  =\Psi\left(  -k_{\ast};W_{0}\right)  =0
\label{S3E3}%
\end{equation}
and with the property that $\Psi\left(  k;W_{0}\right)  \neq0$ for any
$k\in\mathbb{R}$ such that $\left\vert k\right\vert >\left\vert k_{\ast
}\right\vert .$ The asymptotic behaviour of the function $\Psi\left(
k;W_{0}\right)  $ as $\left\vert k\right\vert \rightarrow\infty$ is%
\begin{equation}
\Psi\left(  k;W_{0}\right)  \sim a \ \text{sgn}\left(  k\right)  e^{i\frac
{\pi\text{sgn}\left(  k\right)  }{2}\left(  \frac{\gamma}{2}+p-\frac{1}%
{2}\right)  }\left\vert k\right\vert ^{\frac{1}{2}+\left(  \frac{\gamma}%
{2}+p\right)  }\text{ as }\left\vert k\right\vert \rightarrow\infty
\label{S3E3b}%
\end{equation}
where%
\[
a=\frac{2i}{1+\gamma+2p}\Gamma\left(  \frac{1}{2}-\left(  \frac{\gamma}%
{2}+p\right)  \right)  .
\]

\end{theorem}


\begin{remark}
\label{RealImPart}Theorem \ref{LinearzPb} implies that for $W_{0}$ as in that
Theorem we have
\[
\mathcal{L}\left(  e^{ik_{0}\cdot};W_{0}\right)  \left(  X\right)  \equiv
B\left(  1,e^{ik_{0}\cdot};W_{0}\right)  \left(  X\right)  +B\left(
e^{ik_{0}\cdot},1;W_{0}\right)  \left(  X\right)  =0\ \ \text{for }%
X\in\mathbb{R}.
\]

Given that $W_{0}$ is real in the real line this implies, taking real and
imaginary parts
\begin{align*}
\mathcal{L}\left(  \cos\left(  k_{0}\cdot\right)  ;W_{0}\right)  \left(
X\right)   &  \equiv B\left(  1,\cos\left(  k_{0}\cdot\right)  ;W_{0}\right)
\left(  X\right)  +B\left(  \cos\left(  k_{0}\cdot\right)  ,1;W_{0}\right)
\left(  X\right)  =0\ \ \text{for }X\in\mathbb{R}\\
\mathcal{L}\left(  \sin\left(  k_{0}\cdot\right)  ;W_{0}\right)  \left(
X\right)   &  \equiv B\left(  1,\sin\left(  k_{0}\cdot\right)  ;W_{0}\right)
\left(  X\right)  +B\left(  \sin\left(  k_{0}\cdot\right)  ,1;W_{0}\right)
\left(  X\right)  =0\ \ \text{for }X\in\mathbb{R}.
\end{align*}

\end{remark}

Theorem \ref{LinearzPb} will be proved in Section \ref{sec:proofofTh2.2}.
\medskip

\begin{remark}
It would be possible to prove the results of this paper with functions $W_{0}$
satisfying just some differentiability conditions, instead of the analyticity
condition formulated in Theorem \ref{LinearzPb}. The main reason to use
analytic functions in a wedge extending toward infinity is because this allows
to obtain automatically estimates for the derivatives of $W$. In particular,
the validity of the asymptotic formula (\ref{S3E3b}) in a wedge implies the
validity of asymptotic formulas for the derivatives along the real line.  These
formulas are obtained by just formally differentiating both sides of the asymptotic
formula (\ref{S3E3b}); the results is a consequence of the classical Cauchy
estimates for analytic functions.

Under the analyticity assumption the solutions that we obtain are very regular. In particular, the singularities of the function $\Phi(s)$ as  $s\to 0^+$ or $s\to 1^-$ are given by some power law.  
\end{remark}

\subsection{Bifurcation of non-power law constant flux solutions.}

\label{ssec:bifurcation}

We recall that our goal is to obtain solutions of (\ref{S2E2}) with $B$ as in
(\ref{S2E4}). Due to the invariance of the problem under rescaling of $H$, we can
assume that $b_{\Phi,\gamma}\sqrt{J_{0}}=1.$ Then $H_{s}\left(  X\right)  =1$. 
Furthermore, we choose $k_{\ast}$ as in Theorem \ref{LinearzPb} and then set
$T=\frac{2\pi}{k_{\ast}}.$ Our goal is to construct kernels $W$, close in some
suitable sense to $W_{0}$ (cf.\ the metric (\ref{A1E8})), and nonconstant
functions $H$ having a period $T$ such that
\begin{equation}
B\left(  H,H;W\right)  \left(  X\right)  =J_{0}=\frac{1}{\left(
b_{\Phi,\gamma}\right)  ^{2}}\ \ ,\ \ X\in\mathbb{R} \,.\label{S4E1}%
\end{equation}

Our plan is to obtain $H$ which behaves approximately as $1+s\cos\left(
kX\right)  $ where $s$ is a small real number. Due to the invariance of the
problem under translations, we could equally well obtain functions $H$ behaving
approximately as $1+s\cos\left(  k_{\ast}\left(  X-X_{0}\right)  \right)  $,
with $X_{0}\in\mathbb{R}$ and $s$ small.

We describe now in detail the functional spaces in which the operator $B$
acts. We will assume that the operator acts on spaces of real functions.

We will denote as $H_{per}^{1}\left(  \left[  0,T\right]  \right)  $ the
Hilbert space obtained as the closure of the functions $f:\left[  0,T\right]
\rightarrow\left[  0,T\right]  $ satisfying
\begin{equation}
f\in C^{\infty}\left(  \left[  0,T\right]  \right)  \ \ ,\ \ \ \partial
_{X}^{\ell}f\left(  0\right)  =\partial_{X}^{\ell}f\left(  T\right)
\ \ ,\ \ \text{for any }\ell=0,1,2,3,... \label{A3E3a}%
\end{equation}
with the norm
\begin{equation}
\left\Vert f\right\Vert =\left(  \int_{0}^{T}\left[  \left\vert f\left(
X\right)  \right\vert ^{2}+\left\vert \partial_{X}f\left(  X\right)
\right\vert ^{2}\right]  dX\right)  ^{\frac{1}{2}} \ . \label{A3E3}%
\end{equation}
Notice that we can identify the elements of the space $H_{per}^{1}\left(
\left[  0,T\right]  \right)  $ with the functions $f:\mathbb{R}\rightarrow
\mathbb{R}$ which satisfy $f\left(  X\right)  =f\left(  X+T\right)  $ for each
$X\in\mathbb{R}$ as well as $\left\Vert f\right\Vert <\infty$ with $\left\Vert
\cdot\right\Vert $ as in (\ref{A3E3}). With this identification, we
have
\[
\left\Vert f\right\Vert =\left(  \int_{a}^{a+T}\left[  \left\vert f\left(
X\right)  \right\vert ^{2}+\left\vert \partial_{X}f\left(  X\right)
\right\vert ^{2}\right]  dX\right)  ^{\frac{1}{2}}\text{ \ for each }%
a\in\mathbb{R} .
\]

Morrey's inequality implies that $f$ is H\"{o}lder continuous in the whole
real line as well as the estimate $\left\Vert f\right\Vert _{L^{\infty}\left(
\mathbb{R}\right)  }\leq C\left\Vert f\right\Vert $ for some constant $C$
depending on $T.$

We will use also the spaces $H_{per}^{s}\left(  \left[  0,T\right]  \right)  $
with $s>0.$ These spaces are the closure of (\ref{A3E3a}) with the norm%
\begin{equation}
\left\Vert f\right\Vert _{H_{per}^{s}}=\left(  \sum_{n\in\mathbb{Z}}\left(
1+\left\vert n\right\vert ^{2s}\right)  \left\vert a_{n}\right\vert
^{2}\right)  ^{\frac{1}{2}}\ \ \text{with\ \ }a_{n}=\frac{1}{T}\int_{0}%
^{T}f\left(  X\right)  e^{-\frac{2\pi Xi}{T}}dX. \label{A3E3b}
\end{equation}

Suppose that $k_{\ast}$ is as in Theorem \ref{LinearzPb}. We now introduce two
subspaces of $H_{per}^{1}\left(  \left[  0,T\right]  \right)  $ as follows
\begin{equation}
V_{1}={\text{span}\left\{  \cos\left(  k_{\ast}\cdot\right)
,\sin\left(  k_{\ast}\cdot\right)  \right\}  }\ \ ,\ \ V_{2}=\overline
{\text{span}\left\{  \cos\left(  nk_{\ast}\cdot\right)  ,\sin\left(  nk_{\ast
}\cdot\right)  :n\neq0,-1,1,\ n\in\mathbb{Z}\right\}  } \label{S4E6}%
\end{equation}
where $\text{span}$ denotes finite linear combinations and the closure is with
respect to the topology of $H_{per}^{1}\left(  \left[  0,T\right]  \right)
.$
%
We define subspaces
\begin{align}
&  Z_{0}={\text{span}\left\{  1\right\}  }\ ,\nonumber\\
&  Z_{1}={\text{span}\left\{  \cos\left(  k_{\ast}\cdot\right)  ,\sin\left(
k_{\ast}\cdot\right)  \right\}  }\ \ ,\ \ \nonumber\\
&  Z_{2}=\overline{\text{span}\left\{  \cos\left(  nk_{\ast}\cdot\right)
,\sin\left(  nk_{\ast}\cdot\right)  :n\neq0,-1,1,\ n\in\mathbb{Z}\right\}  }
\label{S4E6a}%
\end{align}
where the closure here is understood with respect to the topology of
$L^{2}\left(  \left[  0,T\right]  \right)  .$ We will denote as $P_{0}%
,P_{1},P_{2}$ the orthogonal projection of $L^{2}\left(  \left[  0,T\right]
\right)  $ into $Z_{0},\ Z_{1},\ Z_{2}$ respectively. These projections are
given by
\[
P_{0}f=\int_{0}^{1}f\left(  sT\right)  ds\ ,\ P_{1}f=2\cos\left(  k_{\ast
}x\right)  \int_{0}^{1}\cos\left(  2\pi s\right)  f\left(  sT\right)
ds+2\sin\left(  k_{\ast}x\right)  \int_{0}^{1}\sin\left(  2\pi s\right)
f\left(  sT\right)  ds
\]%
\[
P_{2}=I-P_{0}-P_{1}.
\]

We now remark the following. Suppose that $W$ is a locally bounded function
that satisfies $W\left(  Y\right)  =W\left(  -Y\right)  $ and (\ref{S2E3}).
Then, using that $\left\Vert f\right\Vert _{L^{\infty}\left(  \mathbb{R}%
\right)  }\leq C\left\Vert f\right\Vert $ we readily obtain that
\[
\left\Vert B\left(  H,H;W\right)  \right\Vert _{L^{2}}\leq C_{1}\left\Vert
B\left(  H,H;W\right)  \right\Vert _{L^{\infty}}\leq C_{2}\left\Vert
H\right\Vert ^{2}%
\]
where $C_{1},C_{2}$ depend on $T$ and $W.$ Our goal is to obtain nonconstant
solutions of (\ref{S4E1}). We now indicate the strategy that we will follow to
prove the existence of such solutions. Notice that due to the homogeneity of
$B\left(  H,H;W\right)  $ in $H,$ if we obtain a solution of the problem
\begin{equation}
B\left(  \tilde{H},\tilde{H};W\right)  =K_{0} \label{S4E2a}%
\end{equation}
for some $K_{0}>0,$ we would then obtain a solution of (\ref{S4E1}) by means
of $H=\frac{\sqrt{J_{0}}}{\sqrt{K_{0}}}\tilde{H}.$ Applying the operators
$P_{1},P_{2}$ into the equation (\ref{S4E2a}) we obtain
\begin{equation}
P_{1}B\left(  \tilde{H},\tilde{H};W\right)  \left(  X\right)  =P_{2}B\left(
\tilde{H},\tilde{H};W\right)  \left(  X\right)  =0. \label{S4E2}
\end{equation}

Our plan is to prove the existence of solutions of (\ref{S4E2}) for some
function $\tilde{H}$ close to $1.$ Using that $B\left(  1,1;W\right)  >0$ it
will then follow using a continuity argument that (\ref{S4E2a}) holds for some
$K_{0}>0.$ More precisely, we look for a solution of the equations
(\ref{S4E2}) with the form
\begin{equation}
\tilde{H}\left(  X\right)  =1+U\left(  X\right)  \ \ ,\ \ U\left(  X\right)
=s\varphi\left(  X\right)  +\psi\left(  X,s\right)  \label{S4E3}%
\end{equation}
where
\begin{equation}
\varphi\in V_{1}\ \ ,\ \ \psi\in V_{2} \label{S4E4}%
\end{equation}
and $\left\vert s\right\vert $ is small and where $W$ is close to $W_{0}$ in
the uniform convergence norm. The solution that we will obtain will satisfy
$\left\Vert \psi\right\Vert \leq Cs^{2}.$ Under these assumptions we will
obtain
\[
P_{0}B\left(  \tilde{H},\tilde{H};W\right)  =P_{0}B\left(  1,1;W\right)
+P_{0}B\left(  1,U;W\right)  +P_{0}B\left(  U,1;W\right)  +P_{0}B\left(
U,U;W\right)
\]
and using that $P_{0}B\left(  1,U;W\right)  =P_{0}B\left(  U,1;W\right)  =0$
as well as the fact that $B\left(  1,1;W\right)  $ is constant, we obtain
\[
P_{0}B\left(  \tilde{H},\tilde{H};W\right)  =B\left(  1,1;W\right)
+P_{0}B\left(  U,U;W\right)  \equiv K_{0} \,.%
\]
Since $W$ is close to $W_{0}$ in the uniform convergence norm, we
obtain that $B\left(  1,1;W\right)  =\frac{1}{\left(  b_{W,\gamma}\right)
^{2}}$ is close to $\frac{1}{\left(  b_{W_{0},\gamma}\right)  ^{2}}>0.$
Therefore, assuming that $\left\vert s\right\vert $ is sufficiently small and
using also (\ref{S4E2}), we obtain that $\tilde{H}$ solves (\ref{S4E2a}) for
some $K_{0}>0.$

We have then reduced the problem (\ref{S4E1}) to finding a solution of
(\ref{S4E2})-(\ref{S4E4}) for $\left\vert s\right\vert $ small and with $\psi$
satisfying $\left\Vert \psi\right\Vert \leq Cs^{2}.$ Due to the invariance of
the problem (\ref{S4E1}) under translations in the variable $X$ we can assume
without loss of generality that $\varphi\left(  X\right)  =\cos\left(
kX\right)  .$

We now introduce an auxiliary linear operator $\mathcal{L}$. For each $W\in
L^{\infty}\left(  \mathbb{R}\right)  $ satisfying $W\left(  Y\right)
=W\left(  -Y\right)  ,\ Y\in\mathbb{R}$ and the asymptotics (\ref{S2E3}) we
define (cf.\ (\ref{S2E9}))
\begin{equation}
\mathcal{L}\left(  \varphi;W\right)  \equiv B\left(  1,\varphi;W\right)
+B\left(  \varphi,1;W\right)  . \label{S4E5}%
\end{equation}

We have that for each $W$ satisfying the previous assumptions the operator
$\mathcal{L}\left(  \cdot;W\right)  $ is well defined in $H_{per}^{1}\left(
\left[  0,T\right]  \right)  .$ Moreover, 
\begin{equation}
\mathcal{L}\left(  V_{1};W\right)  \subset Z_{1}\ ,\ \ \mathcal{L}\left(
V_{2};W\right)  \subset Z_{2} \label{S4E5a}%
\end{equation}
with $V_{1},\ V_{2}$ as in (\ref{S4E6}) and $Z_{1},\ Z_{2}$ are as in
(\ref{S4E6a}). Choosing $W_{0}$ as in Theorem \ref{LinearzPb}, we
have also $\mathcal{L}\left(  V_{1};W_{0}\right)  =\left\{  0\right\}  $ and
in particular
\begin{equation}
\mathcal{L}\left(  \cos\left(  k_{\ast}\cdot\right)  ;W_{0}\right)  =0.
\label{S4E7}
\end{equation}

We now make precise the choice of $W.$ We will choose $W$ with the form
\begin{equation}
W=W_{0}+W_{1} \label{S4E7a}%
\end{equation}
where $W_{1}$ is the linear combination of two functions $W_{1,1},$ $W_{1,2}$
which are analytic in the domain $D_{\beta}\subset\mathbb{C}$, $\beta>0$,
introduced in the statement of Theorem \ref{LinearzPb}
(cf.~\eqref{def:Dbeta}). Moreover, $W_{1,1},$ $W_{1,2}$ tend to zero as
$\left\vert Y\right\vert \rightarrow\infty$ and satisfy
\begin{equation}
W_{1,1}\left(  Y\right)  =W_{1,1}\left(  -Y\right)  \ \ .\ \ W_{1,2}\left(
Y\right)  =W_{1,2}\left(  -Y\right)  \ \ \text{for }Y\in D_{\beta
}.\ \label{S4E8}%
\end{equation}

We have the following result.

\begin{lemma}
\label{FunctW1}Let $k_{\ast}>0$ be as in Theorem \ref{LinearzPb}. There exist
two functions $W_{1,1},\ W_{1,2}$ analytic in $D_{\beta}$ for some $\beta>0,$
real in the real axis, satisfying (\ref{S4E8}) as well as the estimate
\begin{equation}
\left\vert W_{1,1}\left(  Y\right)  \right\vert +\left\vert W_{1,2}\left(
Y\right)  \right\vert \leq\exp\left(  -a\left\vert Y\right\vert ^{2}\right)
\text{ for }Y\in D_{\beta} \label{S4E9}%
\end{equation}
for some $a>0,$ such that
\begin{equation}
\text{span}\left\{  \mathcal{L}\left(  \varphi;W_{1,1}\right)  ,\mathcal{L}%
\left(  \varphi;W_{1,2}\right)  \right\}  =Z_{1} \label{S5E1}%
\end{equation}
with $\varphi_{0}\left(  X\right)  =\cos\left(  k_{\ast}X\right)  .$
\end{lemma}

\begin{proof}
Suppose that $W_{1,1},\ W_{1,2}$ are two functions analytic in $D_{\beta},$
satisfying (\ref{S4E8}), as well as $W_{1,j}\left(  Y\right)  \in\mathbb{R}$
for $Y\in\mathbb{R}$, and that decay sufficiently fast to ensure that
$\left\vert \Psi\left(  k_{\ast};W_{1,1}\right)  \right\vert$, \\ $ \left\vert
\Psi\left(  k_{\ast};W_{1,2}\right)  \right\vert $ defined by means of
(\ref{S3E5}), (\ref{S3E4}) are finite.$\ $Using (\ref{S3E3a}) we can write
\begin{equation}
\mathcal{L}\left(  \varphi;W_{1,j}\right)  =\frac{1}{2}\left[  \Psi\left(
k_{\ast};W_{1,j}\right)  e^{ik_{\ast}\cdot}+\Psi\left(  -k_{\ast}%
;W_{1,j}\right)  e^{-ik_{\ast}\cdot}\right]  \ \ ,\ \ j=1,2\label{S5E2}%
\end{equation}
where $\Psi\left(  k_{\ast};W_{1,j}\right)  $ is as in (\ref{S3E4}),
(\ref{S3E5}). Using that $\Psi\left(  -k_{\ast};W_{1,j}\right)  =\overline
{\Psi\left(  k_{\ast};W_{1,j}\right)  }$ we can rewrite (\ref{S5E2}) as
\[
\mathcal{L}\left(  \varphi;W_{1,j}\right)  =\text{Re}\left(  \Psi\left(
k_{\ast};W_{1,j}\right)  e^{ik_{\ast}\cdot}\right)  \ \ ,\ \ j=1,2.
\]

Writing $\Psi\left(  k_{\ast};W_{1,j}\right)  $ in polar form $\left\vert
\Psi\left(  k_{\ast};W_{1,j}\right)  \right\vert e^{i\arg\left(  \Psi\left(
k_{\ast};W_{1,j}\right)  \right)  }$ we readily see that the functions
$\left\{  \mathcal{L}\left(  \varphi;W_{1,j}\right)  \right\}  _{j=1,2}$ are
linearly independent in $L^{2}\left(  \mathbb{R}\right)  $ iff the vectors
$$\left(
\begin{array}
[c]{c}%
\text{Re}\left(  \Psi\left(  k_{\ast};W_{1,1}\right)  \right) \\
\text{Im}\left(  \Psi\left(  k_{\ast};W_{1,1}\right)  \right)
\end{array}
\right)  ,\ \ \ \left(
\begin{array}
[c]{c}%
\text{Re}\left(  \Psi\left(  k_{\ast};W_{1,2}\right)  \right) \\
\text{Im}\left(  \Psi\left(  k_{\ast};W_{1,2}\right)  \right)
\end{array}
\right)  $$ considered as elements of $\mathbb{R}^{2}$ are linearly
independent. Using (\ref{S3E4}), (\ref{S3E5}) we easily see that in the case
of $W_{1,j}\left(  z\right)  =\delta\left(  z-z_{j}\right)  ,\ j=1,2,$ for
some $z_{1},z_{2}\in\mathbb{R}$ we would have
\[
\Psi\left(  k_{\ast};W_{1,j}\right)  =\frac{G\left(  z_{j},k_{\ast}\right)
}{ik_{\ast}}\ \ ,\ \ j=1,2.
\]

Since the function $z\rightarrow\arg\left(  G\left(  z,k_{\ast}\right)
\right)  $ for a fixed value of $k_{\ast}>0$ is not constant, we obtain the
existence of two values $z_{1},z_{2}$ yielding the desired linear independence
conditions. In order to obtain $W_{1,1},$ $W_{1,2}$ with the desired symmetry,
analyticity conditions and decay at infinity we argue as in the proof of
Theorem \ref{LinearzPb}. More precisely, we replace the Dirac masses
$\delta\left(  z-z_{j}\right)  ,\ j=1,2$ by the functions $\zeta_{\varepsilon
}\left(  z-z_{1}\right)  ,\ \zeta_{\varepsilon}\left(  z-z_{2}\right)  $ where
the function $\zeta_{\varepsilon}$ is as in (\ref{S3E6}) and $\varepsilon>0$.
Using the continuity of the functions $W_{1,j} \mapsto \Psi\left(  k_{\ast};W_{1,j}\right)
,\ j=1,2$, in the weak topology, we
obtain functions $W_{1,1},\ W_{1,2}$ with the properties stated in the Lemma. 
\end{proof}

We can now continue with the analysis of the nonlinear bifurcation problem
which has been reduced to the analysis of (\ref{S4E2})-(\ref{S4E4}) with
$\left\vert s\right\vert $ small. We will prove now the following result.

\begin{theorem}
\label{ProjProb}Suppose that $W_{0},\ k_{\ast}$ are as in Theorem
\ref{LinearzPb} and $W_{1,1},W_{1,2}$ are as in Lemma \ref{FunctW1}. Let
$\varphi_0\left(  X\right)  =\cos\left(  k_{\ast}X\right)  .$ Then, there exists
$s_{0}>0$ and a constant $C>0$ depending only on $W_{0},W_{1,1},W_{1,2}$ such
that, for any $s\in\mathbb{R}$ satisfying $\left\vert s\right\vert \leq s_{0}$
there exist $\alpha_{1},\alpha_{2}\in\mathbb{R}$ and $\psi\in V_{2}$
satisfying $\left\Vert \psi(\cdot,s)\right\Vert \leq C\left\vert s\right\vert
^{2}$ such that the function $\tilde{H}$ defined by means of (\ref{S4E3})
solves (\ref{S4E2}) with
\begin{equation}
W=W_{0}+\alpha_{1}W_{1,1}+\alpha_{2}W_{1,2} \label{S5E5}%
\end{equation}
where%
\begin{equation}
W>0. \label{S5E5a}%
\end{equation}

\end{theorem}

\begin{remark}
The value of $s_{0}$ depends on the function $W_{0}$\ in Theorem
\ref{LinearzPb}. As a matter of fact, the construction of the function $W_{0}$
that will be made in Section \ref{SectionLinear} will depend on several
parameters $z_{a},\ z_{b},\;a,\ b,\ \varepsilon,$ which are assumed to be
fixed. Then, the parameter $s_{0}$ depends on all these parameters.
\end{remark}

\textit{Proof of Theorem \ref{ProjProb}. First part. } We first reformulate the problem
(\ref{S4E2}) as a fixed point problem for a suitable operator. Using
(\ref{S2E4}) as well as the definition of the operator $\mathcal{L}$ in
\eqref{S4E5} we obtain
\[
B\left(  1+U,1+U;W\right)  =B\left(  1,1;W\right)  +\mathcal{L}\left(
U;W\right)  +B\left(  U,U;W\right).
\]

Taking now the operators $P_{1}$ and $P_{2}$ of this expression and using that
$P_{j}B\left(  1,1;W\right)  =0$ for $j=1,2$ we can rewrite (\ref{S4E2}) as
\begin{align*}
P_{1}\mathcal{L}\left(  U;W\right)  +P_{1}B\left(  U,U;W\right)   &  =0\\
P_{2}\mathcal{L}\left(  U;W\right)  +P_{2}B\left(  U,U;W\right)   &  =0.
\end{align*}

Using now (\ref{S4E5a}) and the formula of $U$ in (\ref{S4E3}) we obtain
\begin{equation}
sP_{1}\mathcal{L}\left(  \varphi;W\right)  +P_{1}B\left(  U,U;W\right)
=0\ \ ,\ \ P_{2}\mathcal{L}\left(  \psi;W\right)  +P_{2}B\left(  U,U;W\right)
=0. \label{S6E2}%
\end{equation}

We write $W_{1}=\alpha_{1}W_{1,1}+\alpha_{2}W_{1,2}.$ Then, assuming that $W$
has the form (\ref{S5E5}) we obtain
\[
sP_{1}\mathcal{L}\left(  \varphi_0;W_{0}\right)  +sP_{1}\mathcal{L}\left(
\varphi_0;W_{1}\right)  +P_{1}B\left(  U,U;W\right)  =0.
\]

Using then (\ref{S4E7}) we have $\mathcal{L}\left(  \varphi_0;W_{0}\right)  =0.$
On the other hand, we have that $\mathcal{L}\left(  \varphi_0;W_{1}\right)  \in
Z_{1}.$ Then
\[
s\mathcal{L}\left(  \varphi_0;W_{1}\right)  +P_{1}B\left(  U,U;W\right)  =0
\]
or equivalently
\begin{equation}
\alpha_{1}\mathcal{L}\left(  \varphi_0;W_{1,1}\right)  +\alpha_{2}%
\mathcal{L}\left(  \varphi_0;W_{1,2}\right)  =-\frac{1}{s}P_{1}B\left(
U,U;W_{0}+W_{1}\right).  \label{A1}%
\end{equation}

Due to Lemma \ref{FunctW1} we have that the vectors $\left\{  \mathcal{L}%
\left(  \varphi_0;W_{1,1}\right)  ,\mathcal{L}\left(  \varphi_0;W_{1,2}\right)
\right\}  $ span $Z_{1}.$ Therefore, there exist two linear forms $\ell
_{j}:Z_{1}\rightarrow\mathbb{R}$, $j=1,2,$ such that for any $w\in Z_{1}$ we
have
\[
w=\ell_{1}\left(  w\right)  \mathcal{L}\left(  \varphi_0;W_{1,1}\right)
+\ell_{2}\left(  w\right)  \mathcal{L}\left(  \varphi_0;W_{1,2}\right).
\]

Notice that $\left\vert \ell_{1}\left(  w\right)  \right\vert +\left\vert
\ell_{2}\left(  w\right)  \right\vert \leq C\left\Vert w\right\Vert
_{L^{1}\left(  \mathbb{R}\right)  }.$ The constant $C$ depends on the
functions $W_{1,1},$ $W_{1,2}$ but these will be assumed to be assumed to be
fixed in all the remaining argument.

Then, due to (\ref{A1}), it follows that the constants $\alpha_{1},\alpha_{2}$
in (\ref{S5E5}) must be chosen as
\begin{equation}
\alpha_{1}=-\frac{1}{s}\ell_{1}\left(  P_{1}B\left(  U,U;W_{0}+W_{1}\right)
\right)  \ \ ,\ \ \alpha_{2}=-\frac{1}{s}\ell_{2}\left(  P_{1}B\left(
U,U;W_{0}+W_{1}\right)  \right).  \label{S6E1}%
\end{equation}

This choice of $\alpha_{1},\alpha_{2}$ implies the first equation in
(\ref{S6E2}). Notice that (\ref{S6E1}) is an equation for the coefficients
$\alpha_{1},\alpha_{2}$ because $W_{1}$ depends on both of them. Nevertheless,
it has a suitable form for a fixed point argument. It only remains to
reformulate the second equation in (\ref{S6E2}). This requires to examine the
invertibility properties of the operator $P_{2}\mathcal{L}$.

The choice of $k_{\ast}$ (cf.\ Theorem \ref{LinearzPb}) implies that
$\Psi\left(  nk_{\ast};W\right)  \neq0$ for any $n\neq0,\ \pm1.$ Indeed, this
follows from the fact that $\Psi\left(  k;W_{0}\right)  \neq0$ for $\left\vert
k\right\vert >k_{\ast}$ as well as the fact that the functions $\Psi\left(
k;W\right)  $ depend continuously on the function $W$ in the uniform topology
of measures. Moreover, the asymptotics of $\Psi\left(  k;W_{0}\right)  $ as
$\left\vert k\right\vert \rightarrow\infty$ is given by (\ref{S3E3b}) and this
asymptotic behaviour is not modified adding to $W_{0}$ the function
$\alpha_{1}W_{1,1}+\alpha_{2}W_{1,2}.$ Therefore, the claim follows if $s_{0}$
is small enough.

The operator $P_{2}\mathcal{L}\left(  \cdot;W\right)  $ acts in Fourier as
follows. A function $f\in Z_{2}$ can be represented by means of a Fourier
series with the form
\[
f\left(  X\right)  =\sum_{n\in\mathbb{Z}\setminus\left\{  0,\pm1\right\}
}a_{n}e^{ink_{\ast}X}\ \ \text{with }\sum_{n\in\mathbb{Z}\setminus\left\{
0,\pm1\right\}  }n^{2}\left\vert a_{n}\right\vert ^{2}<\infty
\]
where in addition we have $\overline{a_{n}}=a_{-n}.$ Then, using (\ref{S3E3a})
we would obtain
\begin{equation}
A_{W}\left(  f\right)  \left(  X\right)  :=P_{2}\mathcal{L}\left(  f;W\right)
\left(  X\right)  =\sum_{n\in\mathbb{Z}\setminus\left\{  0,\pm1\right\}  }%
\Psi\left(  nk_{\ast};W\right)  a_{n}e^{ink_{\ast}X}. \label{A2}%
\end{equation}

The expression (\ref{A2}) defines an operator acting on $f\in H_{per}%
^{1}\left(  \left[  0,T\right]  \right)  \cap Z_{2}$ for each function $W$ as
in (\ref{S5E5}).

Using (\ref{S3E3b}) it then follows that $A_{W}\left(  f\right)  \in
H_{per}^{\frac{1}{2}-\left(  \frac{\gamma}{2}+p\right)  }\left(  \left[
0,T\right]  \right)  \cap\left\{  a_{0}=a_{1}=a_{-1}=0\right\}  =H_{per}%
^{\frac{1}{2}-\left(  \frac{\gamma}{2}+p\right)  }\left(  \left[  0,T\right]
\right)  \cap Z_{2}$ where the spaces $H_{per}^{\frac{1}{2}-\left(
\frac{\gamma}{2}+p\right)  }$ are endowed with the norm (\ref{A3E3b}).
Therefore, the inverse $\left(  A_{W}\right)  ^{-1}$ is defined in the space
$H_{per}^{\frac{1}{2}-\left(  \frac{\gamma}{2}+p\right)  }\left(  \left[
0,T\right]  \right)  \cap Z_{2}$ and it transforms this space in a subspace of
$H_{per}^{1}\left(  \left[  0,T\right]  \right)  .$ Moreover, we have%
\begin{equation}
\left\Vert \left(  A_{W}\right)  ^{-1}\left(  \varphi\right)  \right\Vert
_{H_{per}^{1}\left(  \left[  0,T\right]  \right)  }\leq C\left\Vert
\varphi\right\Vert _{H_{per}^{\frac{1}{2}-\left(  \frac{\gamma}{2}+p\right)
}\left(  \left[  0,T\right]  \right)  } \label{A4E8}%
\end{equation}
for each $\varphi\in H_{per}^{\frac{1}{2}-\left(  \frac{\gamma}{2}+p\right)
}\left(  \left[  0,T\right]  \right)  \cap Z_{2}$.

We can now rewrite the second equation in (\ref{S6E2}) in a more convenient
form in order to reformulate (\ref{S6E2}) as a fixed point problem for
$\left(  \alpha_{1},\alpha_{2},\psi\right)  $. More precisely, suppose that we
have $B\left(  U,U;W\right)  \in H_{per}^{\frac{1}{2}-\left(  \frac{\gamma}%
{2}+p\right)  }\left(  \left[  0,T\right]  \right)  $ for any $U\in
H_{per}^{1}\left(  \left[  0,T\right]  \right)  .$ We can then write
\[
P_{2}\mathcal{L}\left(  \psi;W_{0}\right)  +P_{2}\mathcal{L}\left(  \psi
;W_{1}\right)  =-P_{2}B\left(  U,U;W\right)
\]%
\begin{equation}
\psi=-\left(  A_{W}\right)  ^{-1}\left(  P_{2}\mathcal{L}\left(  \psi
;W_{1}\right) \right)  -\left(  A_W\right)  ^{-1}\left(
P_{2}B\left(  U,U;W_{0}+W_{1}\right)  \right).  \label{S6E3}%
\end{equation}

Notice that (\ref{S6E1}), (\ref{S6E3}) yield a reformulation of the problem
(\ref{S6E2}) as a fixed point problem for $\left(  \alpha_{1},\alpha_{2}%
,\psi\right)  $, assuming that $W_{1}=\alpha_{1}W_{1,1}+\alpha_{2}%
W_{1,2},\ U=s\varphi+\psi.$ More precisely, if we define the mapping%
\begin{equation}
T\left(
\begin{array}
[c]{c}%
\alpha_{1}\\
\alpha_{2}\\
\psi
\end{array}
\right)  =\left(
\begin{array}
[c]{c}%
-\frac{1}{s}\ell_{1}\left(  P_{1}B\left(  U,U;W_{0}+W_{1}\right)  \right)  \\
-\frac{1}{s}\ell_{2}\left(  P_{1}B\left(  U,U;W_{0}+W_{1}\right)  \right)  \\
-\left(  A_{W}\right)  ^{-1}\left(  P_{2}\mathcal{L}\left(  \psi;W_{1}\right)
\right)  -\left( A_W\right)  ^{-1}\left(  P_{2}B\left(
U,U;W_{0}+W_{1}\right) \right)
\end{array}
\right)  \label{S6E3b}%
\end{equation}
we have that (\ref{S6E1}), (\ref{S6E3}) can be rewritten as%
\begin{equation}
\left(
\begin{array}
[c]{c}%
\alpha_{1}\\
\alpha_{2}\\
\psi
\end{array}
\right)  =T\left(
\begin{array}
[c]{c}%
\alpha_{1}\\
\alpha_{2}\\
\psi
\end{array}
\right).  \label{S6E3a}%
\end{equation}


In order to conclude the proof of Theorem \ref{ProjProb} we must check that
the operator in (\ref{S6E3b}) is well defined for $U\in H_{per}^{1}\left(
\left[  0,T\right]  \right)  .$ Specifically we need to prove that $B\left(
U,U;W\right)  $ transforms $U\in H_{per}^{1}\left(  \left[  0,T\right]
\right)  $ into $H_{per}^{\frac{1}{2}-\left(  \frac{\gamma}{2}+p\right)
}\left(  \left[  0,T\right]  \right)  .$ To this end we prove the following Lemma.

\begin{lemma}
\label{bilinEst}Let $k_{\ast}$ be as in Theorem \ref{LinearzPb}. Suppose that
$W$ is analytic in a domain $D_{\beta}$ for some $\beta>0.$ Suppose that
(\ref{S3E2}) holds, with $\gamma+2p<1.$ Let $B$ be as in (\ref{S2E4}). Then,
the following estimate holds
\begin{equation}
\left\Vert B\left(  U_{1},U_{2};W\right)  \right\Vert _{H_{per}^{\frac{1}%
{2}-\left(  \frac{\gamma}{2}+p\right)  }\left(  \left[  0,T\right]  \right)
}\leq C\left\Vert U_{1}\right\Vert \left\Vert U_{2}\right\Vert \ \ ,\ U_{1}%
,U_{2}\in H_{per}^{1}\left(  \left[  0,T\right]  \right)  \label{S6E5}%
\end{equation}
where $C$ depends only on $\gamma,p,\kappa,k_{\ast}$, the multiplicative
constant in the term $O\left(  e^{-\kappa\left\vert X\right\vert }\right)  $
and the supremum of $W$ in $D_{\beta}.$
\end{lemma}

\begin{remark}
Notice that $W_{0}$ satisfies also (\ref{S3E2}) due to \eqref{S5E5} and Lemma
\ref{FunctW1}.
\end{remark}

\begin{proof}
Suppose that $U_{1},U_{2}\in H_{per}^{1}\left(  \left[  0,T\right]  \right)
.$ We can represent them as
\[
U_{1}\left(  X\right)  =\sum_{n\in\mathbb{Z}}a_{n}e^{inkX}\ ,\ \ U_{2}\left(
X\right)  =\sum_{n\in\mathbb{Z}}b_{n}e^{inkX}\ \ ,\ \ \overline{a_{n}}%
=a_{-n}\ \ ,\ \ \overline{b_{n}}=b_{-n}\ ,\
\]
where we recall that $k_{\ast}=\frac{2\pi}{T}.$ Notice that the series
defining $U_{1}$ converges absolutely since $\sum_{n\in\mathbb{Z}}\left\vert
a_{n}\right\vert \leq\left(  \sum_{n\in\mathbb{Z}}n^{2}\left\vert
a_{n}\right\vert ^{2}\right)  ^{\frac{1}{2}}\left(  \sum_{n\in\mathbb{Z}}%
\frac{1}{n^{2}}\right)  ^{\frac{1}{2}}\leq C\left\Vert U_{1}\right\Vert .$ The
same argument applies to $U_{2}.$

Then, using the definition of the operator $B$ we obtain%
\begin{align}
&  B\left(  U_{1},U_{2};W\right)  \left(  X\right) \label{S6E6}\\
&  =\sum_{n\in\mathbb{Z}}\sum_{m\in\mathbb{Z}}a_{m}b_{n}\int_{-\infty}%
^{X}dY\int_{X+\log\left(  1-e^{Y-X}\right)  }^{\infty}dZ\left[  e^{\frac{1}%
{2}\left(  Y-Z\right)  }W\left(  Y-Z\right)  \right]  e^{imk_{\ast}%
Y}e^{ink_{\ast}Z}\nonumber\\
&  =\sum_{n\in\mathbb{Z}}\sum_{m\in\mathbb{Z}}J_{m,n}a_{m}b_{n}e^{i\left(
m+n\right)  k_{\ast}X}=\sum_{\ell\in\mathbb{Z}}\left[  \sum_{m+n=\ell}%
J_{m,n}a_{m}b_{n}\right]  e^{i\ell k_{\ast}X}\nonumber
\end{align}
where
\[
J_{m,n}=\int_{-\infty}^{0}dY\int_{\log\left(  1-e^{Y}\right)  }^{\infty
}dZ\left[  e^{\frac{1}{2}\left(  Y-Z\right)  }W\left(  Y-Z\right)  \right]
e^{imk_{\ast}Y}e^{ink_{\ast}Z}.
\]
Notice that the coefficients $J_{m,n}$ are well defined due to (\ref{S3E2}).
We can compute them, using the change of variables $\xi=Z-Y,\ d\xi=dZ$ and
applying Fubini's Theorem as follows
\begin{align}
J_{m,n}  &  =\int_{-\infty}^{0}e^{imk_{\ast}Y}dY\int_{\log\left(
1-e^{Y}\right)  }^{\infty}dZ\left[  e^{-\frac{1}{2}\left(  Z-Y\right)
}W\left(  Z-Y\right)  \right]  e^{ink_{\ast}Z}\nonumber\\
&  =\int_{-\infty}^{0}e^{i\left(  m+n\right)  k_{\ast}Y}dY\int_{-Y+\log\left(
1-e^{Y}\right)  }^{\infty}\left[  e^{-\frac{\xi}{2}}e^{ink_{\ast}\xi}W\left(
\xi\right)  \right]  d\xi\nonumber\\
&  =\int_{-\infty}^{\infty}e^{-\frac{\xi}{2}}e^{ink_{\ast}\xi}W\left(
\xi\right)  d\xi\int_{-\log\left(  e^{\xi}+1\right)  }^{0}e^{i\left(
m+n\right)  k_{\ast}Y}dY\nonumber\\
&  =\frac{1}{i\left(  m+n\right)  k_{\ast}}\int_{-\infty}^{\infty}%
e^{-\frac{\xi}{2}}e^{ink_{\ast}\xi}W\left(  \xi\right)  \left[  1-\frac
{1}{\left(  e^{\xi}+1\right)  ^{i\left(  m+n\right)  k_{\ast}}}\right]
d\xi\label{S6E7}%
\end{align}
if $m+n\neq0$, and
\begin{equation}
J_{m,n}=\int_{-\infty}^{\infty}e^{-\frac{\xi}{2}}e^{ink_{\ast}\xi}W\left(
\xi\right)  \log\left(  e^{\xi}+1\right)  d\xi\text{ \ \ if }m+n=0.
\label{S6E7a}%
\end{equation}

Due to (\ref{S6E6}), in order to estimate the Fourier coefficients of
$B\left(  U_{1},U_{2};W\right)  $ we need to derive bounds for the sums
\begin{equation}
\sum_{m+n=\ell}J_{m,n}a_{m}b_{n}. \label{S6E8}%
\end{equation}

Using (\ref{S6E7}) we obtain that the coefficients $J_{m,n}$ with
$m+n=\ell\neq0$ are given by
\begin{equation}
J_{m,n}=\frac{1}{i\ell k_{\ast}}\int_{-\infty}^{\infty}e^{-\frac{\xi}{2}%
}e^{ink_{\ast}\xi}W\left(  \xi\right)  \left[  1-\frac{1}{\left(  e^{\xi
}+1\right)  ^{i\ell k_{\ast}}}\right]  d\xi\ \ ,\ \ m+n=\ell.\label{S6E9}%
\end{equation}

In order to estimate these coefficients we write
\begin{equation}
J_{m,n}=J_{m,n}^{\left(  1\right)  }+J_{m,n}^{\left(  2\right)  } \label{S7E1}%
\end{equation}
where
\begin{align*}
J_{m,n}^{\left(  1\right)  }  &  =\frac{1}{i\ell k_{\ast}}\int_{-\infty}%
^{0}e^{-\frac{\xi}{2}}e^{ink_{\ast}\xi}W\left(  \xi\right)  \left[  1-\frac
{1}{\left(  e^{\xi}+1\right)  ^{i\ell k_{\ast}}}\right]  d\xi\\
J_{m,n}^{\left(  2\right)  }  &  =\frac{1}{i\ell k_{\ast}}\int_{0}^{\infty
}e^{-\frac{\xi}{2}}e^{ink_{\ast}\xi}W\left(  \xi\right)  \left[  1-\frac
{1}{\left(  e^{\xi}+1\right)  ^{i\ell k_{\ast}}}\right]  d\xi.
\end{align*}

Using (\ref{S3E2}) and that $k_{\ast}>0$, we obtain
\begin{equation}
\left\vert J_{m,n}^{\left(  2\right)  }\right\vert \leq\frac{C}{\left\vert
\ell\right\vert k_{\ast}}\ \ ,\ \ \ell\neq0. \label{S7E2}%
\end{equation}

In order to estimate $J_{m,n}^{\left(  1\right)  }$ we use the change of
variables $\eta=\log\left(  1+e^{\xi}\right)  ,$ whence $\xi=\log\left(
e^{\eta}-1\right)  $ and $d\xi=\frac{e^{\eta}}{e^{\eta}-1}d\eta.$ Then
\begin{align*}
J_{m,n}^{\left(  1\right)  }  &  =\frac{1}{i\ell k_{\ast}}\int_{0}%
^{\log\left(  2\right)  }\frac{W\left(  \log\left(  e^{\eta}-1\right)
\right)  }{\left(  e^{\eta}-1\right)  ^{\frac{1}{2}}}\left(  e^{\eta
}-1\right)  ^{ink_{\ast}}\left[  1-e^{-i\ell k_{\ast}\eta}\right]
\frac{e^{\eta}}{e^{\eta}-1}d\eta\\
&  =\frac{1}{i\ell k_{\ast}}\int_{0}^{\log\left(  2\right)  }\frac{e^{\eta
}W\left(  \log\left(  e^{\eta}-1\right)  \right)  }{\left(  e^{\eta}-1\right)
^{\frac{3}{2}-ink_{\ast}}}\left(  1-e^{-i\ell k_{\ast}\eta}\right)  d\eta.
\end{align*}

The analyticity properties of $W$ imply that the function $Q\left(
\eta\right)  =e^{\eta}W\left(  \log\left(  e^{\eta}-1\right)  \right)  $ is
analytic in a domain containing the interval $\left(  0,\log\left(  2\right)
\right]  .$ The domain contains a whole neighbourhood of the point $\eta
=\log\left(  2\right)  .$ Concerning the neighbourhood of $\eta=0$ we can see,
using Taylor series for $e^{\eta}$ that we have analyticity of $Q$ if
$\left\vert \eta\right\vert $ is small enough and $\arg\log\left(
\eta\right)  \in\left(  -\pi-\beta,-\pi+\beta\right)  $ with $\beta>0,$ i.e.
the region of analyticity of $Q$ covers the whole set $\left\vert
\eta\right\vert \leq\delta$ with $\delta$ small, $\arg\left(  \delta\right)
\in\left(  -\pi,\pi\right)  .$ Moreover, the function $Q$ can be extended
analytically along the negative real axis in a small neighbourhood of the
origin to yield a multivalued function.

The asymptotic behaviour of $Q$ near the origin can be obtained using
(\ref{S3E2}). We have%
\[
Q\left(  \eta\right)  =\left(  \eta\right)  ^{-\left(  \frac{\gamma}%
{2}+p\right)  }\left[  1+O\left(  \left\vert \eta\right\vert ^{\kappa}\right)
\right]  \text{ as }\left\vert \eta\right\vert \rightarrow0.
\]

Therefore%
\[
\left\vert Q\left(  \eta\right)  \right\vert \leq C\left\vert \eta\right\vert
^{-\left(  \frac{\gamma}{2}+p\right)  }\text{ for }0<\eta\leq\log\left(
2\right).
\]
We then write
\[
J_{m,n}^{\left(  1\right)  }=\frac{1}{i\ell k_{\ast}}\int_{0}^{\log\left(
2\right)  }\frac{Q\left(  \eta\right)  }{\left(  e^{\eta}-1\right)  ^{\frac
{3}{2}-ink_{\ast}}}\left(  1-e^{-i\ell k_{\ast}\eta}\right)  d\eta=\frac
{1}{i\ell k_{\ast}}\int_{0}^{\log\left(  2\right)  }\frac{Q\left(
\eta\right)  }{\left(  e^{\eta}-1\right)  ^{\frac{3}{2}}}\frac{\left(
1-e^{-i\ell k_{\ast}\eta}\right)  }{\left(  e^{\eta}-1\right)  ^{-ink_{\ast}}%
}d\eta.
\]

Using that $\left\vert \left(  e^{\eta}-1\right)  ^{-ink_{\ast}}\right\vert
=1$ and $\left\vert 1-e^{-i\ell k_{\ast}\eta}\right\vert \leq C\left\vert
\ell\right\vert \eta$ for $\eta\in\left(  0,\log\left(  2\right)  \right)  $
and $\ell\in\mathbb{Z}$ we obtain%
\[
\left\vert J_{m,n}^{\left(  1\right)  }\right\vert \leq\frac{C}{\left\vert
\ell\right\vert }\int_{0}^{\log\left(  2\right)  }\frac{\left\vert Q\left(
\eta\right)  \right\vert }{\left(  e^{\eta}-1\right)  ^{\frac{3}{2}}}%
\frac{\left\vert \ell\right\vert \eta}{\left\vert \left(  e^{\eta}-1\right)
^{-ink_{\ast}}\right\vert }d\eta\leq C\int_{0}^{\log\left(  2\right)  }%
\frac{d\eta}{\left\vert \eta\right\vert ^{\left(  \frac{\gamma}{2}+p+\frac
{1}{2}\right)  }}\leq C.
\]

Combining this with (\ref{S7E2}) we obtain
\[
\left\vert J_{m,n}\right\vert \leq C\ \ \text{for }m+n=\ell\neq0.
\]

On the other hand, (\ref{S6E7a}) combined with (\ref{S7E2}) implies
\[
\left\vert J_{m,n}\right\vert \leq C\ \ \text{for }m+n=0\ .
\]
Therefore
\begin{equation}
\left\vert J_{m,n}\right\vert \leq C\ \label{S7E4}%
\end{equation}
for $m$ and $n$ arbitrary integers.

We can now estimate the $\ell-$th Fourier coefficient of $B\left(  U_{1}%
,U_{2};W\right)  $ which is given in (\ref{S6E8}) as
\begin{equation}
\left\vert \sum_{m+n=\ell}J_{m,n}a_{m}b_{n}\right\vert \leq C\sum_{m+n=\ell
}\left\vert a_{m}\right\vert \left\vert b_{n}\right\vert \ . \label{S7E5}%
\end{equation}
We have
\[
\left\Vert B\left(  U_{1},U_{2};W\right)  \right\Vert _{H_{per}^{\frac{1}%
{2}-\left(  \frac{\gamma}{2}+p\right)  }\left(  \mathbb{R}\right)  }^{2}%
=\sum_{\ell\in\mathbb{Z}}\left[  \left(  1+\left\vert \ell\right\vert \right)
^{\left(  1-\gamma-2p\right)  }\left\vert \sum_{m+n=\ell}J_{m,n}a_{m}%
b_{n}\right\vert ^{2}\right]
\]
and using (\ref{S7E5}) we obtain
\begin{equation}
\left\Vert B\left(  U_{1},U_{2};W\right)  \right\Vert _{H_{per}^{\frac{1}%
{2}-\left(  \frac{\gamma}{2}+p\right)  }\left(  \mathbb{R}\right)  }^{2}\leq
C\sum_{\ell\in\mathbb{Z}}\left(  1+\left\vert \ell\right\vert \right)
^{2}\left(  \sum_{n\in\mathbb{Z}}\left\vert a_{\ell-n}\right\vert \left\vert
b_{n}\right\vert \right)  ^{2}. \label{S7E6}%
\end{equation}

In order to estimate the right-hand side of (\ref{S7E6}) we define the
following periodic functions in $S^{1}:$%
\[
\psi_{1}\left(  x\right)  =\sum_{n\in\mathbb{Z}}\left\vert a_{n}\right\vert
e^{inx}\ \ ,\ \ \psi_{2}\left(  x\right)  =\sum_{n\in\mathbb{Z}}\left\vert
b_{n}\right\vert e^{inx}\ \ ,\ \ x\in\left[  -\pi,\pi\right]  \ .
\]
Using Plancherel formula as well as the fact that $U_{1},U_{2}$ are in
$H_{per}^{1}\left(  \mathbb{R}\right)  $ we can estimate the right-hand side
of (\ref{S7E6}) as
\begin{equation}
C\left[  \left\Vert \psi_{1}\psi_{2}\right\Vert _{L^{2}\left(  -\pi
,\pi\right)  }^{2}+\left\Vert \partial_{x}\left(  \psi_{1}\psi_{2}\right)
\right\Vert _{L^{2}\left(  -\pi,\pi\right)  }^{2}\right]  \ . \label{S7E7}%
\end{equation}
Moreover,
\begin{align*}
\left\Vert \psi_{1}\right\Vert _{L^{2}\left(  -\pi,\pi\right)  }%
^{2}+\left\Vert \partial_{x}\left(  \psi_{1}\right)  \right\Vert
_{L^{2}\left(  -\pi,\pi\right)  }^{2}  &  \leq C\left\Vert U_{1}\right\Vert
^{2}\\
\left\Vert \psi_{2}\right\Vert _{L^{2}\left(  -\pi,\pi\right)  }%
^{2}+\left\Vert \partial_{x}\left(  \psi_{2}\right)  \right\Vert
_{L^{2}\left(  -\pi,\pi\right)  }^{2}  &  \leq C\left\Vert U_{2}\right\Vert
^{2}%
\end{align*}
and
\[
\left\Vert \psi_{j}\right\Vert _{L^{\infty}\left(  -\pi,\pi\right)  }\leq
C\left\Vert U_{j}\right\Vert \text{ \ , \ }j=1,2\ .
\]

Then, the term in (\ref{S7E7}) can be estimated by $C\left\Vert U_{1}%
\right\Vert ^{2}\left\Vert U_{2}\right\Vert ^{2}.$ Therefore, using also
(\ref{S7E6}) we obtain
\[
\left\Vert B\left(  U_{1},U_{2};W\right)  \right\Vert _{H_{per}^{\frac{1}%
{2}-\left(  \frac{\gamma}{2}+p\right)  }\left(  \mathbb{R}\right)  }^{2}\leq
C\left\Vert U_{1}\right\Vert ^{2}\left\Vert U_{2}\right\Vert ^{2}%
\]
whence the Lemma follows.
\end{proof}

We can estimate easily the dependence of the bilinear operator in $W_{1}.$

\begin{lemma}
\label{bilincont}Let $\alpha\in\mathbb{R}^{2}$ with $\alpha=\left(  \alpha
_{1},\alpha_{2}\right)  $ and $W_{1}=\alpha_{1}W_{1,1}+\alpha_{2}W_{1,2}$ with
$W_{1,1},\ W_{1,2}$ as in Lemma \ref{FunctW1}. Then, the following estimate
holds%
\begin{equation}
\left\Vert B\left(  U_{1},U_{2};W_{1}\right)  \right\Vert _{H_{per}^{\frac
{1}{2}-\left(  \frac{\gamma}{2}+p\right)  }\left(  \left[  0,T\right]
\right)  }\leq C\left\vert \alpha\right\vert \left\Vert U_{1}\right\Vert
\left\Vert U_{2}\right\Vert \ \ ,\ U_{1},U_{2}\in H_{per}^{1}\left(  \left[
0,T\right]  \right).  \label{A4E9}%
\end{equation}

\end{lemma}

\begin{proof}
Suppose that $U_{1},U_{2}\in H_{per}^{1}\left(  \left[  0,T\right]  \right)
.$ We represent then using a Fourier series as in the Proof of Lemma
\ref{bilinEst}
\[
U_{1}\left(  X\right)  =\sum_{n\in\mathbb{Z}}a_{n}e^{ink_{\ast}X}%
\ ,\ \ U_{2}\left(  X\right)  =\sum_{n\in\mathbb{Z}}b_{n}e^{ink_{\ast}%
X}\ \ ,\ \ \overline{a_{n}} =a_{-n}\ \ ,\ \ \overline{b_{n}}=b_{-n}.
\]

Then, arguing as in the Proof of Lemma \ref{bilinEst} we obtain
\begin{equation}
B\left(  U_{1},U_{2};W_{1}\right)  \left(  X\right)  =\sum_{\ell\in\mathbb{Z}%
}\left[  \sum_{m+n=\ell}J_{m,n}a_{m}b_{n}\right]  e^{i\ell k_{\ast}X}\nonumber
\end{equation}
where $J_{m,n}$ is given in \eqref{S6E7}-\eqref{S6E7a}.

We have $\left\vert J_{m,n}\right\vert \leq C$ if $\ell=m+n=0.$ To estimate
$J_{m,n}$ for $\ell\neq0$ we write%
\begin{equation}
J_{m,n}=J_{m,n}^{\left(  1\right)  }+J_{m,n}^{\left(  2\right)  }%
\end{equation}
where
\begin{align*}
J_{m,n}^{\left(  1\right)  }  &  =\frac{1}{i\ell k_{\ast}}\int_{-\infty}%
^{0}e^{-\frac{\xi}{2}}e^{ink_{\ast}\xi}W_{1}\left(  \xi\right)  \left[
1-\frac{1}{\left(  e^{\xi}+1\right)  ^{i\ell k_{\ast}}}\right]  d\xi\\
J_{m,n}^{\left(  2\right)  }  &  =\frac{1}{i\ell k_{\ast}}\int_{0}^{\infty
}e^{-\frac{\xi}{2}}e^{ink_{\ast}\xi}W_{1}\left(  \xi\right)  \left[
1-\frac{1}{\left(  e^{\xi}+1\right)  ^{i\ell k_{\ast}}}\right]  d\xi.
\end{align*}

Arguing as in the Proof of Lemma \ref{bilinEst} we obtain%
\[
\left\vert J_{m,n}^{\left(  2\right)  }\right\vert \leq\frac{C}{\left\vert
\ell\right\vert k_{\ast}}\left\vert \alpha\right\vert \ \ ,\ \ \ell\neq0.
\]

Moreover, due to the fast decay of $W_{1}\left(  \xi\right)  $ as
$\xi\rightarrow-\infty$ (cf.\ (\ref{S4E9})) we can obtain a similar estimate
for $J_{m,n}^{\left(  1\right)  }$. Therefore%
\[
\left\vert J_{m,n}^{\left(  1\right)  }\right\vert \leq\frac{C}{\left\vert
\ell\right\vert k_{\ast}}\left\vert \alpha\right\vert \ \ ,\ \ \ell\neq0
\]
whence, combining all the estimates we obtain%
\begin{equation}
\left\vert J_{m,n}\right\vert \leq\frac{C}{1+\left\vert \ell\right\vert
}\left\vert \alpha\right\vert \ \ ,\ \ \ell\neq0\ \ \text{for }m+n=\ell
\label{A5E1}%
\end{equation}
where $C$ depends on $k_{\ast}$ but not on $m,n.$ We then obtain the estimate
arguing as in the Proof of Lemma \ref{bilinEst}. Indeed, the proof of that
Lemma just relies on the boundedness of $\left\vert J_{m,n}\right\vert .$ The
estimate (\ref{A5E1}) implies that $\left\vert J_{m,n}\right\vert \leq
C\left\vert \alpha\right\vert .$ Using this, a simple adaptation of the
argument in the Proof of Lemma \ref{bilinEst} yields (\ref{A4E9}) whence the
result follows.
\end{proof}

Lemma \ref{bilincont} yields also estimates for the linear operator
$\mathcal{L}\left(  U;W_{1}\right)  .$

\begin{lemma}
\label{linW1}Let $k_{\ast}$ be as in Theorem \ref{LinearzPb} and $\alpha,$
$W_{1}$ as in Lemma \ref{bilinEst}. Let $\mathcal{L}\left(  \varphi
;W_{1}\right)  $ be as in (\ref{S4E5}) with $W=W_{1}$ with $\varphi\in
H_{per}^{1}\left(  \left[  0,T\right]  \right)  .$ Then%
\[
\left\Vert \mathcal{L}\left(  \varphi;W_{1}\right)  \right\Vert _{H_{per}%
^{\frac{1}{2}-\left(  \frac{\gamma}{2}+p\right)  }\left(  \left[  0,T\right]
\right)  }\leq C\left\vert \alpha\right\vert \left\Vert \varphi\right\Vert
\ \ ,\ \varphi\in H_{per}^{1}\left(  \left[  0,T\right]  \right)
\]
where $C$ depends on $k_{\ast}$ but it is independent of $\alpha$ and
$\varphi.$
\end{lemma}

\begin{proof}
It is just a consequence of the fact that $\mathcal{L}\left(  \varphi
;W_{1}\right)  \equiv B\left(  1,\varphi;W_{1}\right)  +B\left(
\varphi,1;W_{1}\right)  $ and (\ref{A4E9}) in Lemma \ref{bilincont}.
\end{proof}

We can conclude now the Proof of Theorem \ref{ProjProb}.

\begin{proofof}
[End of the Proof of Theorem \ref{ProjProb}]The problem can be reformulated as
a fixed point. To this end we introduce a Hilbert space
\[
\mathcal{H}=\left\{  \left(  \alpha_{1},\alpha_{2},\psi\right)  :\alpha
_{1},\alpha_{2}\in\mathbb{R}\text{,\ }\left\Vert \psi\right\Vert
<\infty\right\}
\]
as well as a mapping $T:\mathcal{H}\rightarrow\mathcal{H}$%
\begin{equation}
T\left(  \alpha_{1},\alpha_{2},\psi\right)  =\left(  T_{1},T_{2},T_{3}\right)
\left(  \alpha_{1},\alpha_{2},\psi\right)  \in\mathcal{H}\label{A5E2}%
\end{equation}
where
\begin{align}
T_{1}\left(  \alpha_{1},\alpha_{2},\psi\right)   &  =-\frac{1}{s}\ell
_{1}\left(  P_{1}B\left(  U,U;W_{0}+W_{1}\right)  \right)  \label{A5E3}\\
T_{2}\left(  \alpha_{1},\alpha_{2},\psi\right)   &  =-\frac{1}{s}\ell
_{2}\left(  P_{1}B\left(  U,U;W_{0}+W_{1}\right)  \right)  \nonumber\\
T_{3}\left(  \alpha_{1},\alpha_{2},\psi\right)   &  =-\left(  A_{W}\right)
^{-1}\left(  P_{2}\mathcal{L}\left(  \psi;W_{1}\right)  \right)
-\left(  A_{W}\right)  ^{-1}\left(  P_{2}B\left(  U,U;W_{0}+W_{1}\right)
\right)  \nonumber
\end{align}
with
\[
U=s\varphi_{0}+\psi\ \ ,\ \ W_{1}=\alpha_{1}W_{1,1}+\alpha_{2}W_{1,2}.
\]

Then, the problem (\ref{S6E1}), (\ref{S6E3}) is equivalent to the fixed point
problem
\begin{equation}
\label{eq:defT}\left(  \alpha_{1},\alpha_{2},\psi\right)  =T\left(  \alpha
_{1},\alpha_{2},\psi\right).
\end{equation}

We now define the set $K_{M,s_{0}}=\left\{  \left(  \alpha_{1},\alpha_{2}%
,\psi\right)  :\left\vert \alpha_{1}\right\vert +\left\vert \alpha
_{2}\right\vert \leq M\left\vert s\right\vert ,\ \left\Vert \psi\right\Vert
\leq\left\vert s\right\vert ,\ \left\vert s\right\vert \leq s_{0}\right\}  .$
We introduce in $K_{M,s_{0}}$ the metric
\begin{equation}
\mathrm{dist}\left(  \left(  \alpha_{1}^{\left(  1\right)  },\alpha
_{2}^{\left(  1\right)  },\psi^{\left(  1\right)  }\right)  ,\left(
\alpha_{1}^{\left(  2\right)  },\alpha_{2}^{\left(  2\right)  },\psi^{\left(
2\right)  }\right)  \right)  =\left\vert \alpha_{1}^{\left(  1\right)
}-\alpha_{1}^{\left(  2\right)  }\right\vert +\left\vert \alpha_{2}^{\left(
1\right)  }-\alpha_{2}^{\left(  2\right)  }\right\vert +M\left\Vert
\psi^{\left(  1\right)  }-\psi^{\left(  2\right)  }\right\Vert. \label{dist}%
\end{equation}

We will show that if $M$ is chosen sufficiently large and $s_{0}$ sufficiently
small the operator $T$ defined by means of (\ref{A5E2}), (\ref{A5E3}) is
contractive. We emphasize here that the whole argument is made assuming that
the function $W_{0}$ in (\ref{S5E5}) is fixed. The function $W_{0}$ depends on
$\varepsilon$ and then all the constants $C$ in the following might depend on
$\varepsilon.$ However, we will choose the constants $C$ independent of $M$
and $s.$ In the following argument, we need to assume that $M$ is sufficiently
large (depending on $\varepsilon$) and then $s_{0}$ sufficiently small
(depending on $M$).

We first prove that $T$ transforms $K_{M,s_{0}}$ into itself. To this end we
derive estimates for $T_{1}\left(  \alpha_{1},\alpha_{2},\psi\right)
,\ T_{2}\left(  \alpha_{1},\alpha_{2},\psi\right)  $ for each $\left(
\alpha_{1},\alpha_{2},\psi\right)  \in K_{M,s_{0}}.$ Assuming that
$M\left\vert s\right\vert $ is sufficiently small we can apply Lemma
\ref{bilinEst} to estimate $\left\Vert B\left(  U,U;W_{0}+W_{1}\right)
\right\Vert _{H_{per}^{\frac{1}{2}-\left(  \frac{\gamma}{2}+p\right)  }\left(
\mathbb{R}\right)  }$ as $C\left\vert s\right\vert ^{2}$ with $C$ independent
of $M.$ Then $\left\vert T_{1}\left(  \alpha_{1},\alpha_{2},\psi\right)
\right\vert +\left\vert T_{2}\left(  \alpha_{1},\alpha_{2},\psi\right)
\right\vert \leq C_{0}\left\vert s\right\vert $ with $C_{0}$ independent of
$M.$ Therefore, choosing $M>C_{0}$ we obtain that the two components
$T_{1}\left(  \alpha_{1},\alpha_{2},\psi\right)  ,\ T_{2}\left(  \alpha
_{1},\alpha_{2},\psi\right)  $ satisfy the two inequalities required in the
definition of $K_{M,s_{0}}.$ On the other hand, using Lemma \ref{linW1} and
the fact that $P_{2}$ is a projection operator, we can estimate $P_{2}%
\mathcal{L}\left(  \psi;W_{1}\right)  $ as $\left\Vert P_{2}\mathcal{L}\left(
\psi;W_{1}\right)  \right\Vert _{H_{per}^{\frac{1}{2}-\left(  \frac{\gamma}%
{2}+p\right)  }\left(  \left[  0,T\right]  \right)  }\leq CM\left\vert
s\right\vert ^{2}$ with $C$ independent of $M.$ On the other hand, Lemma
\ref{bilinEst} implies $\left\Vert B\left(  U,U;W_{0}+W_{1}\right)
\right\Vert _{H_{per}^{\frac{1}{2}-\left(  \frac{\gamma}{2}+p\right)  }\left(
\left[  0,T\right]  \right)  }\leq C\left\vert s\right\vert ^{2}$ assuming
that $M\left\vert s\right\vert $ is sufficiently small. Using (\ref{A4E8}) we
obtain%
\begin{align}
\left\Vert \left(  A_{W}\right)  ^{-1}\left(  P_{2}\mathcal{L}\left(
\psi;W_{1}\right)  \right)  \right\Vert  &  \leq CM\left\vert
s\right\vert ^{2}\label{eq:P2L_s2}\\
\left\Vert \left(  A_{W}\right)  ^{-1}\left(  P_{2}B\left(  U,U;W_{0}%
+W_{1}\right)  \right)  \right\Vert  &  \leq C\left\vert s\right\vert
^{2}.\label{eq:P2L_s2_2}%
\end{align}

Therefore, if we choose $\left\vert s\right\vert $ small enough we obtain
$\left\Vert T_{3}\left(  \alpha_{1},\alpha_{2},\psi\right)  \right\Vert
\leq\left\vert s\right\vert .$ This estimate combined with the estimates for
$T_{1}\left(  \alpha_{1},\alpha_{2},\psi\right)  ,\ T_{2}\left(  \alpha
_{1},\alpha_{2},\psi\right)  $ obtained above imply that for each $\left(
\alpha_{1},\alpha_{2},\psi\right)  \in K_{M,s_{0}}$ we have $T\left(
\alpha_{1},\alpha_{2},\psi\right)  \in K_{M,s_{0}}.$

It remains to show that the operator $T$ is contractive. Given $\left(
\alpha_{1}^{\left(  j\right)  },\alpha_{2}^{\left(  j\right)  },\psi^{\left(
j\right)  }\right)  \in K_{M,s_{0}},\ j=1,2,$ we write%
\[
U^{\left(  j\right)  }=s\varphi+\psi^{\left(  j\right)  }\ \ ,\ \ W_{1}%
=\alpha_{1}^{\left(  j\right)  }W_{1,1}+\alpha_{2}^{\left(  j\right)  }%
W_{1,2}\ ,\ \alpha^{\left(  j\right)  }=\left(  \alpha_{1}^{\left(  j\right)
},\alpha_{2}^{\left(  j\right)  }\right)  \ ,\ \ \ j=1,2.
\]

Then, using that $B\left(  U,U;W\right)  $ is multilinear in its arguments,
and using also Lemmas \ref{bilinEst}, \ref{bilincont} we obtain%
\begin{align*}
&  \left\Vert B\left(  U^{\left(  1\right)  },U^{\left(  1\right)  }%
;W_{0}+W_{1}^{\left(  1\right)  }\right)  -B\left(  U^{\left(  2\right)
},U^{\left(  2\right)  };W_{0}+W_{1}^{\left(  2\right)  }\right)  \right\Vert
_{H_{per}^{\frac{1}{2}-\left(  \frac{\gamma}{2}+p\right)  }\left(
\mathbb{R}\right)  }\\
&  \leq C\left\vert s\right\vert \left\Vert \psi^{\left(  1\right)  }%
-\psi^{\left(  2\right)  }\right\Vert +C\left\vert s\right\vert ^{2}\left\vert
\alpha^{\left(  1\right)  }-\alpha^{\left(  2\right)  }\right\vert
\end{align*}
where $C$ is independent of $M.$ Then, using the first two equations of
(\ref{A5E3}) as well as (\ref{dist}) we obtain%
\begin{align}
&  \left\vert T_{k}\left(  \alpha_{1}^{\left(  1\right)  },\alpha_{2}^{\left(
1\right)  },\psi^{\left(  1\right)  }\right)  -T_{k}\left(  \alpha
_{1}^{\left(  2\right)  },\alpha_{2}^{\left(  2\right)  },\psi^{\left(
2\right)  }\right)  \right\vert \nonumber\\
&  \leq C\left\Vert \psi^{\left(  1\right)  }-\psi^{\left(  2\right)
}\right\Vert +C\left\vert s\right\vert \left\vert \alpha^{\left(  1\right)
}-\alpha^{\left(  2\right)  }\right\vert \nonumber \\
& \leq C\left(  \frac{1}{M}+\left\vert
s\right\vert \right)  \mathrm{dist}\left(  \left(  \alpha_{1}^{\left(
1\right)  },\alpha_{2}^{\left(  1\right)  },\psi^{\left(  1\right)  }\right)
,\left(  \alpha_{1}^{\left(  2\right)  },\alpha_{2}^{\left(  2\right)  }
,\psi^{\left(  2\right)  }\right)  \right) \label{A5E4}
\end{align}
for $k=1,2.$ On the other hand, using the last equation in (\ref{A5E3})
combined with Lemmas \ref{bilinEst}, \ref{bilincont} and \ref{linW1} we obtain%
\begin{align}
&  \left\Vert T_{3}\left(  \alpha_{1}^{\left(  1\right)  },\alpha_{2}^{\left(
1\right)  },\psi^{\left(  1\right)  }\right)  -T_{3}\left(  \alpha
_{1}^{\left(  2\right)  },\alpha_{2}^{\left(  2\right)  },\psi^{\left(
2\right)  }\right)  \right\Vert \label{A5E5}\\
&  \leq CM\left\vert s\right\vert \left\Vert \psi^{\left(  1\right)  }%
-\psi^{\left(  2\right)  }\right\Vert +C\left\vert s\right\vert \left\vert
\alpha^{\left(  1\right)  }-\alpha^{\left(  2\right)  }\right\vert
+C\left\vert s\right\vert \left\Vert \psi^{\left(  1\right)  }-\psi^{\left(
2\right)  }\right\Vert +C\left\vert s\right\vert ^{2}\left\vert \alpha
^{\left(  1\right)  }-\alpha^{\left(  2\right)  }\right\vert. \nonumber
\end{align}

Multiplying (\ref{A5E5}) by $M,$ adding the result to (\ref{A5E4}) and using
(\ref{dist}) we obtain%
\begin{align*}
&  \mathrm{dist}\left(  T\left(  \alpha_{1}^{\left(  1\right)  },\alpha
_{2}^{\left(  1\right)  },\psi^{\left(  1\right)  }\right)  ,T\left(
\alpha_{1}^{\left(  2\right)  },\alpha_{2}^{\left(  2\right)  },\psi^{\left(
2\right)  }\right)  \right) \\
&  \leq C\left(  \frac{1}{M}+M\left\vert s\right\vert \right)  \mathrm{dist}%
\left(  \left(  \alpha_{1}^{\left(  1\right)  },\alpha_{2}^{\left(  1\right)
},\psi^{\left(  1\right)  }\right)  ,\left(  \alpha_{1}^{\left(  2\right)
},\alpha_{2}^{\left(  2\right)  },\psi^{\left(  2\right)  }\right)  \right).
\end{align*}

Then, choosing $M$ sufficiently large and then $\left\vert s\right\vert \leq
s_{0}$ we obtain that the operator $T$ is contractive, whence the result
follows. Notice that if $s_{0}$ is sufficiently small we have that $W>0.$ The
quadratic estimate of $\psi$ follows from the estimates
\eqref{eq:P2L_s2}-\eqref{eq:P2L_s2_2} and from the definition of $T$ \eqref{eq:defT}.
\end{proofof}

\medskip

\begin{remark}
\label{ConstFluxSol}Notice that Theorem \ref{ProjProb} implies the existence
of nonconstant solutions of (\ref{S4E1}) for each $J_{0}>0$ and kernels
$W=W_{0}+W_{1}$ with $\left\vert s\right\vert \leq s_{0}.$ Indeed, we have
already seen that the value of $J_{0}$ can be assumed to be any positive
number by means of a rescaling argument. Then, given that $B\left(
1,1;W_{0}\right)  =a>0,$ we obtain that the solutions $\tilde{H}$ of
(\ref{S4E2}) obtained in Theorem \ref{ProjProb} satisfy $B\left(  \tilde
{H},\tilde{H};W\right)  =J_{0}$ with $J_{0}>0$ if $s_{0}$ is small enough due
to the continuous dependence of $B\left(  \tilde{H},\tilde{H};W_{0}\right)  $
on $s.$
\end{remark}

We can now reformulate Theorem \ref{ProjProb} in terms of the original set of
variables (cf.\ (\ref{S1E1}), (\ref{A1E3})) in order to prove Theorem
\ref{main}.

\begin{proofof}
[Proof of Theorem \ref{main}]It is just a consequence of Theorem
\ref{ProjProb} and Remark \ref{ConstFluxSol} using the change of variables
(\ref{S2E1}). Notice that the function $X\rightarrow e^{\frac{\left(
\gamma+3\right)  X}{2}}f\left(  e^{X}\right)  =\tilde{H}\left(  X\right)  $ is
periodic with period $T=\frac{2\pi}{k_{\ast}},$ with $k_{\ast}$ as in Theorem
\ref{LinearzPb}, whence%
\[
e^{\frac{\left(  \gamma+3\right)  X}{2}}f\left(  e^{X}\right)  =e^{\frac
{\left(  \gamma+3\right)  T}{2}}e^{\frac{\left(  \gamma+3\right)  X}{2}%
}f\left(  e^{T}e^{X}\right)  \ .
\]
Therefore
\[
f\left(  Qx\right)  =\frac{f\left(  x\right)  }{Q^{\frac{\left(
\gamma+3\right)  }{2}}}\ \ \ \text{for each } \ x>0 \ \text{ with\ }\ Q=e^{ T}
\ .
\]
We have $f\left(  x\right)  =\tilde{H}\left(  \log\left(  x\right)  \right)
x^{-\frac{\left(  \gamma+3\right)  }{2}}.$ Notice that the function $f$ is not
a power law, since $\tilde{H}$ is not constant for $s\neq0.$ Finally, the
continuity of the family of kernels $K$ in the topology induced by the metric
\eqref{A1E8} follows from \eqref{S2E3a}, the asymptotics of $W_{0}$ in
\eqref{S3E2}, the estimate \eqref{S4E9} and equation \eqref{S5E5} as well as
the continuity in $s$ of the functions $\alpha_{1}$, $\alpha_{2}$ in
\eqref{S5E5}. Then the result follows.
\end{proofof}


\section{Construction of the bifurcation kernel: proof of Theorem
\ref{LinearzPb}\label{SectionLinear}}

\label{sec:proofofTh2.2} In order to prove Theorem \ref{LinearzPb} we first
need an auxiliary result which allows to study the properties of some
complex-valued functions.

\medskip

\begin{lemma}
\label{ZeroesLimitFunct} We define a function $G:\mathbb{R}\times
\mathbb{R}\rightarrow\mathbb{C}$, $G=G\left(  z,k\right)  $,
by means of
\begin{equation}
G\left(  z,k\right)  =e^{-\frac{z}{2}}\left(  1+e^{ikz}\right)  \left(
1-\frac{1}{\left(  e^{z}+1\right)  ^{ik}}\right)  +e^{\frac{z}{2}}\left(
1+e^{-ikz}\right)  \left(  1-\frac{1}{\left(  e^{-z}+1\right)  ^{ik}}\right)
. \label{S3E5}%
\end{equation}
There exists $\delta_{0}>0$ sufficiently small such that, for any $\delta
\in\left(  0,\delta_{0}\right)  $ there exists $z_{0}=z_{0}\left(
\delta\right)  >0$ sufficiently large such that, if $z_{b}\geq z_{0}$ and
$z_{a}=\left(  1+\delta\right)  z_{b}$ there exist infinitely many values
$k_{n}\in\mathbb{R}$, with $n\in\mathbb{N}$ such that%
\begin{equation}
\frac{G\left(  z_{a},k_{n}\right)  }{G\left(  z_{b},k_{n}\right)  }%
\in\mathbb{R} \ \ \text{\ with\ \ }\frac{G\left(  z_{a},k_{n}\right)
}{G\left(  z_{b},k_{n}\right)  }<0 . \label{S3E5a}%
\end{equation}
Moreover, for each $\delta$ and $z_{b}$ there exist a number $\sigma>0$ and
sequences of positive numbers $\left\{  \varepsilon_{n,1}\right\}
_{n\in\mathbb{N}},\ \left\{  \varepsilon_{n,2}\right\}  _{n\in\mathbb{N}}$
such that%
\begin{align}
\arg\left(  G\left(  z_{b},k_{n}-\varepsilon_{n,1}\right)  \right)
-\arg\left(  G\left(  z_{a},k_{n}-\varepsilon_{n,1}\right)  \right)   &
=\pi-\sigma\label{S3E5b}\\
\arg\left(  G\left(  z_{b},k_{n}+\varepsilon_{n,2}\right)  \right)
-\arg\left(  G\left(  z_{a},k_{n}+\varepsilon_{n,2}\right)  \right)   &
=\pi+\sigma\label{S3E5c}%
\end{align}
and also%
\begin{align}
0  &  <R_{1}e^{\frac{z_{a}}{2}}\leq\left\vert G\left(  z_{a},k\right)
\right\vert \leq R_{2}e^{\frac{z_{a}}{2}}<\infty\ \ ,\ \ 0<R_{1}e^{\frac
{z_{a}}{2}}\leq\left\vert G\left(  z_{a},k\right)  \right\vert \leq
R_{2}e^{\frac{z_{a}}{2}}<\infty\label{S3E5d}\\
0  &  <R_{1}e^{\frac{z_{b}}{2}}\leq\left\vert G\left(  z_{b},k\right)
\right\vert \leq R_{2}e^{\frac{z_{b}}{2}}<\infty\ \ ,\ \ 0<R_{1}e^{\frac
{z_{b}}{2}}\leq\left\vert G\left(  z_{b},k\right)  \right\vert \leq
R_{2}e^{\frac{z_{b}}{2}}<\infty\label{S3E5e}%
\end{align}
for $k\in\left[  k_{n}-\varepsilon_{n,1},k_{n}+\varepsilon_{n,2}\right]  $ for
some $R_{1},\ R_{2}$ independent of $z_{a},\ z_{b}.$
\end{lemma}

\begin{proofof}
[Proof of Lemma \ref{ZeroesLimitFunct}]We define $G_{1}\left(  z,k\right)
=e^{\frac{z}{2}}\left(  1+e^{-ikz}\right)  \left(  1-\frac{1}{\left(
e^{-z}+1\right)  ^{ik}}\right)  .$ Using that $\left\vert e^{ikz}\right\vert
=\left\vert \left(  e^{z}+1\right)  ^{ik}\right\vert =1$ we obtain%
\begin{equation}
\left\vert G\left(  z,k\right)  -G_{1}\left(  z,k\right)  \right\vert
\leq4e^{-\frac{z_{0}}{2}}\ \ \text{for all }\left(  z,k\right)  \in
\mathbb{R}\times\mathbb{R}\ \ ,\ \ z\geq z_{0}. \label{A2E0}%
\end{equation}

The function $G_{1}\left(  z,k\right)  $ can be written as%
\begin{equation}
G_{1}\left(  z,k\right)  =e^{\frac{z}{2}}\left(  1+e^{-ikz}\right)  \left(
1-\exp\left(  -ik\omega\left(  z\right)  \right)  \right)  \label{A2E1}%
\end{equation}
where%
\begin{equation}
\omega\left(  z\right)  =\log\left(  1+e^{-z}\right).  \label{A2E2}%
\end{equation}

Notice that (\ref{A2E1}) implies that $G_{1}\left(  z,k\right)  =e^{\frac
{z}{2}}\zeta_{1}\zeta_{2}$ where $\zeta_{1},\zeta_{2}\in\left\{  \zeta
\in\mathbb{C}:\left\vert \zeta-1\right\vert =1\right\}  .$ As $k$ varies from
$-\infty$ to $+\infty$ we have that $\zeta_{1}$ and $\zeta_{2}$ rotate in a
circle of radius $1$ centered at the point $\zeta=1$ in the counterclock sense
with angular velocity $z$ and $\omega\left(  z\right)  $ respectively. We will
denote
\begin{align*}
\zeta_{1,a}=\left(  1+e^{-ikz_{a}}\right)  ,\ \ \ \zeta_{1,b}= \left(
1+e^{-ikz_{b}}\right)  , \ \ \ \zeta_{2,a}=\left(  1+e^{-ik\omega\left(
z_{a}\right)  }\right)  , \ \ \ \zeta_{2,b}=\left(  1+e^{-ik\omega\left(
z_{b}\right)  }\right)  \ .
\end{align*}
Notice that all points $\left(  \zeta_{1,a}-1\right)  ,\ \left(  \zeta
_{1,b}-1\right)  ,\ \left(  \zeta_{2,a}-1\right)  ,\ \left(  \zeta
_{2,b}-1\right)  ,\ \frac{\left(  \zeta_{1,a}-1\right)  }{\left(  \zeta
_{1,b}-1\right)  },\ \frac{\left(  \zeta_{2,b}-1\right)  }{\left(  \zeta
_{2,a}-1\right)  }$ are in the unit circle $\left\{  \zeta\in\mathbb{C}%
:\left\vert \zeta\right\vert =1\right\}  $ for any $k\in\mathbb{R}$. The six
points rotate around the origin in the clockwise sense as $k$ increases with
constant angular velocity $z_{a},\ z_{b},\ \omega\left(  z_{a}\right)
,\ \omega\left(  z_{b}\right)  ,\ \left(  z_{a}-z_{b}\right)  =\delta
z_{b},\ \left(  \omega\left(  z_{b}\right)  -\omega\left(  z_{a}\right)
\right)  $ respectively. This implies that the six points $\left(  \zeta
_{1,a}-1\right)  ,\ \left(  \zeta_{1,b}-1\right)  $ , $\left(  \zeta
_{2,a}-1\right)  ,\ \left(  \zeta_{2,b}-1\right)  ,\ \frac{\left(  \zeta
_{1,a}-1\right)  }{\left(  \zeta_{1,b}-1\right)  },\ \frac{\left(  \zeta
_{2,b}-1\right)  }{\left(  \zeta_{2,a}-1\right)  }$ reach arbitrary positions
in the unit circle in periods $\frac{2\pi}{z_{a}},\ \frac{2\pi}{z_{b}}%
,\ \frac{2\pi}{\omega\left(  z_{a}\right)  },\ \frac{2\pi}{\omega\left(
z_{b}\right)  },\ \frac{2\pi}{\delta z_{b}},\ \frac{2\pi}{\left(
\omega\left(  z_{b}\right)  -\omega\left(  z_{a}\right)  \right)  }$
respectively. In particular choosing $\delta$ sufficiently small and $z_{0}$
sufficiently large (depending on $\delta$) we would obtain, using (\ref{A2E2})
that%
\begin{equation}
\frac{2\pi}{z_{a}}\approx\frac{2\pi}{z_{b}}\ll\frac{2\pi}{\delta z_{b}}%
\ll\frac{2\pi}{\left(  \omega\left(  z_{b}\right)  -\omega\left(
z_{a}\right)  \right)  }\approx\frac{2\pi}{\omega\left(  z_{b}\right)  }%
\ll\frac{2\pi}{\omega\left(  z_{a}\right)  } \label{A2E3} 
\end{equation}
for $z\in (z_b, (1+\delta)z_b)$ with $z_b \geq z_0$ and $\delta$ sufficiently small and $z_{0}$ sufficiently large. 
The inequalities (\ref{A2E3}) imply that for sufficiently large values of $k$
the points $\left(  \zeta_{1,a}-1\right)  ,\ \left(  \zeta_{1,b}-1\right)  $
make an arbitrary angle. On the other hand they rotate with a much faster
angular speed $z_{a}\approx z_{b}$. Therefore, for any prescribed angle
$\alpha\in\left(  0,2\pi\right)  $ and any $\varepsilon_{0}>0$, if we choose
$\delta$ sufficiently small and $z_{0}$ large there are infinitely many values
$k$ for which the points $\left(  \zeta_{1,a}-1\right)  $ can be at any
position of the unit circle with the angle between $\left(  \zeta
_{1,a}-1\right)  $ and$\ \left(  \zeta_{1,b}-1\right)  $ is contained in
$\left(  \alpha-\varepsilon_{0},\alpha+\varepsilon_{0}\right)  .$ The
estimates (\ref{A2E3}) imply also that for $k$ sufficiently large we can make
a similar claim about the points $\left(  \zeta_{2,a}-1\right)  ,\ \left(
\zeta_{2,b}-1\right)  .$ Moreover, (\ref{A2E3}) has also the consequence that
for small $\delta$ and large $z_{0}$ the points $\left(  \zeta_{2,a}-1\right)
,\ \left(  \zeta_{2,b}-1\right)  $ rotate much more slowly as $k$ increases as
the points $\left(  \zeta_{1,a}-1\right)  ,\ \left(  \zeta_{1,b}-1\right)  $
and $\frac{\left(  \zeta_{1,a}-1\right)  }{\left(  \zeta_{1,b}-1\right)  }.$

Therefore, for any $\alpha,\beta,\gamma\in\left(  0,2\pi\right)  $ and any
$\varepsilon_{0}>0$ we can choose $\delta$ sufficiently small and $z_{0}$
large such there are infinitely many values of $k$ such that $\left(
\zeta_{2,a}-1\right)  \in\left(  \gamma-\varepsilon_{0},\gamma+\varepsilon
_{0}\right)  ,\ \left(  \zeta_{1,a}-1\right)  $ is at any position of the
circle, the angle between $\left(  \zeta_{1,a}-1\right)  $ and$\ \left(
\zeta_{1,b}-1\right)  $ is contained in the interval $\left(  \alpha
-\varepsilon_{0},\alpha+\varepsilon_{0}\right)  $ and the angle between
$\left(  \zeta_{2,a}-1\right)  $ and $\left(  \zeta_{2,b}-1\right)  $ is in
$\left(  \beta-\varepsilon_{0},\beta+\varepsilon_{0}\right)  .$

We denote as $\xi_{1},\ \xi_{2}$ two points of the circle $\left\{  \zeta
\in\mathbb{C}:\left\vert \zeta-1\right\vert =1\right\}  $ such that $\xi
_{1}=\rho_{1}e^{i\theta_{1}},\ \xi_{2}=\rho_{2}e^{i\theta_{2}}$ with
$0<\theta_{1}<\frac{\pi}{2},\ -\frac{\pi}{2}<\theta_{2}<0$ and such that
$2\theta_{1}-2\theta_{2}=\pi-s$ where $s>0$ is a small number. Notice that
$\rho_{1}>0,\ \rho_{2}>0.$ Notice that $\xi_{1}=1+e^{i\varphi_{1}},\ \xi
_{2}=1+e^{i\varphi_{2}}$ with $0<\varphi_{1}<\pi$ and $-\pi<\varphi_{2}<0.$ We
have $\theta_{j}=\frac{\varphi_{j}}{2},\ j=1,2$ as it can be seen noticing
that the vector connecting the points $0$ and $1+\exp(i\varphi_{j})$ is a
diagonal of the parallelogram consisting of the points $\{0,1,\exp
(i\varphi_{j}),1+\exp(i\varphi_{j})\}$. Notice that $\rho_{1}$ and $\rho_{2}$
are independent of $z_{a},\ z_{b}$ and $k.$


We choose $\varepsilon_{0}>0$ small and we select $k=k_{0}$ such that
$\zeta_{1,b}=\xi_{1},$ $\left\vert \arg\left(  \zeta_{2,b}-\xi_{1}\right)
\right\vert <\varepsilon_{0},$\ $\left\vert \arg\left(  \zeta_{1,a}-\xi
_{2}\right)  \right\vert <\varepsilon_{0},$ $\left\vert \arg\left(
\zeta_{2,a}-\xi_{2}\right)  \right\vert <\varepsilon_{0}.$ Notice that we need
to take $b$ sufficiently small and $z_{0}$ sufficiently large to ensure that
infinitely many of such a values $k_{0}$ exist and that these values diverge
to $+\infty.$ Therefore, using that $G_{1}\left(  z,k\right)  =e^{\frac{z}{2}%
}\zeta_{1}\zeta_{2}$ it follows that, assuming that $\varepsilon_{0}$ has been
chosen sufficiently small, we have that $\arg\left(  G_{1}\left(  z_{b}%
,k_{0}\right)  \right)  -\arg\left(  G_{1}\left(  z_{a},k_{0}\right)  \right)
$
\[
\pi-2s<\arg\left(  G_{1}\left(  z_{b},k_{0}\right)  \right)  -\arg\left(
G_{1}\left(  z_{a},k_{0}\right)  \right)  <\pi-\frac{s}{2} .
\]

Using (\ref{A2E0}) it then follows that, if $z_{0}$ is chosen sufficiently
large (depending only on $s$,\ $\theta_{1},\ \theta_{2}$) we have%
\begin{equation}
\pi-3s<\arg\left(  G\left(  z_{b},k_{0}\right)  \right)  -\arg\left(  G\left(
z_{a},k_{0}\right)  \right)  <\pi-\frac{s}{4}. \label{A2E4}%
\end{equation}

We now examine the values of $G\left(  z_{a},k\right)  ,\ G\left(
z_{b},k\right)  $ for $k>k_{0}.$ Notice that as $k$ increases we obtain that
$\frac{\zeta_{1,b}}{\zeta_{1,a}}=\frac{e^{-ikz_{b}}}{e^{-ikz_{a}}}=e^{ik\delta
z_{b}}.$ Then, in times of order $\frac{s}{\delta z_{b}}$ we obtain that
$\arg\left(  \zeta_{1,b}\right)  -\arg\left(  \zeta_{1,b}\right)  $ increases
an amount of order $s.$ Then, using the fact that $\zeta_{2,a},\ \zeta_{2,b}$
change much more slowly than $\frac{\zeta_{1,b}}{\zeta_{1,a}}$ (cf.
(\ref{A2E3})), we obtain that there exists $k_{1}\in\left(  k_{0},k_{0}%
+\frac{5s}{\delta z_{b}}\right)  $ such that%
\[
\pi+\frac{s}{2}<\arg\left(  G_{1}\left(  z_{b},k_{1}\right)  \right)
-\arg\left(  G_{1}\left(  z_{a},k_{1}\right)  \right)  <\pi+2s
\]
and in addition $\frac{1}{2}\left\vert G_{1}\left(  z_{b},k\right)
\right\vert \leq\left\vert G_{1}\left(  z_{b},k_{0}\right)  \right\vert
\leq2\left\vert G_{1}\left(  z_{b},k\right)  \right\vert $ and $\frac{1}%
{2}\left\vert G_{1}\left(  z_{a},k\right)  \right\vert \leq\left\vert
G_{1}\left(  z_{a},k_{0}\right)  \right\vert \leq2\left\vert G_{1}\left(
z_{a},k\right)  \right\vert $ for all $k\in\left[  k_{0},k_{0}+k_{1}\right]
.$ Then, using (\ref{A2E0}) it follows that if $z_{0}$ is sufficiently large
we have%
\begin{equation}
\pi+\frac{s}{4}<\arg\left(  G\left(  z_{b},k_{1}\right)  \right)  -\arg\left(
G\left(  z_{a},k_{1}\right)  \right)  <\pi+3s \label{A2E5}%
\end{equation}
and also%
\begin{equation}
\frac{1}{3}\left\vert G\left(  z_{b},k\right)  \right\vert \leq\left\vert
G\left(  z_{b},k_{0}\right)  \right\vert \leq3\left\vert G\left(
z_{b},k\right)  \right\vert \ \ \ \text{and\ \ }\frac{1}{3}\left\vert G\left(
z_{a},k\right)  \right\vert \leq\left\vert G\left(  z_{a},k_{0}\right)
\right\vert \leq3\left\vert G\left(  z_{a},k\right)  \right\vert \label{A2E6}%
\end{equation}
for all $k\in\left[  k_{0},k_{0}+k_{1}\right]  .$ Moreover, since
$G_{1}\left(  z,k\right)  =e^{\frac{z}{2}}\zeta_{1}\zeta_{2},$ $\zeta_{1,a}$
and $\zeta_{2,a}$\ are close to $\xi_{2}$\ and $\zeta_{1,b}$ and $\zeta_{2,b}%
$\ are close to $\xi_{1}$ and $\left\vert \xi_{1}\right\vert $ and $\left\vert
\xi_{2}\right\vert $ are strictly positive (and independent of $k,\ z_{a}%
,\ z_{b}$), it then follows that $\left\vert G_{1}\left(  z_{a},k\right)
\right\vert \geq C_{0}e^{\frac{z_{a}}{2}},\ \left\vert G_{1}\left(
z_{b},k\right)  \right\vert \geq C_{0}e^{\frac{z_{b}}{2}}$ if $k\in\left[
k_{0},k_{0}+k_{1}\right]  ,$ for some $C_{0}>0$ (independent of $k,\ z_{a}%
,\ z_{b}$). Then, combining (\ref{A2E4}) and (\ref{A2E5}) it follows that, if
$z_{0}$ is chosen sufficiently large, there exists at least one value
$k_{\ast}\in\left(  k_{0},k_{0}+k_{1}\right)  $ such that $\arg\left(
G\left(  z_{b},k_{\ast}\right)  \right)  =\arg\left(  G\left(  z_{a},k_{\ast
}\right)  \right)  +\pi.$ The existence of the sequences $\left\{
\varepsilon_{n,1}\right\}  _{n\in\mathbb{N}},\ \left\{  \varepsilon
_{n,2}\right\}  _{n\in\mathbb{N}}$ satisfying (\ref{S3E5b}), (\ref{S3E5c})
with, say $\sigma=\frac{s}{3}$ follows from a similar argument. Estimates
(\ref{S3E5d}), (\ref{S3E5e}) follow from (\ref{A2E0}), (\ref{A2E1}) and the
way in which $\zeta_{1,a},\ \zeta_{1,b},\ \zeta_{2,a},\ \zeta_{2,b}$ have been
chosen. This concludes the proof of the Lemma.
\end{proofof}

\begin{figure}[ptb]
\centering
\includegraphics[scale=0.5]{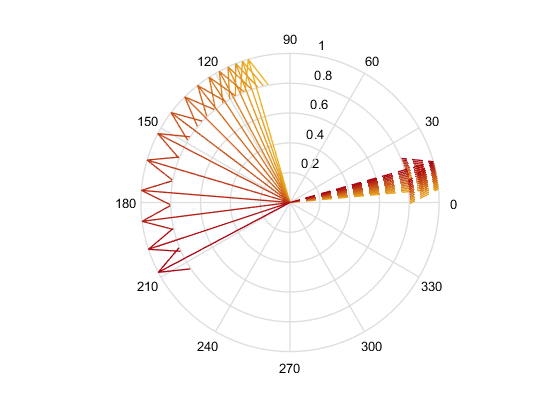}\ \ \includegraphics[scale=0.5]{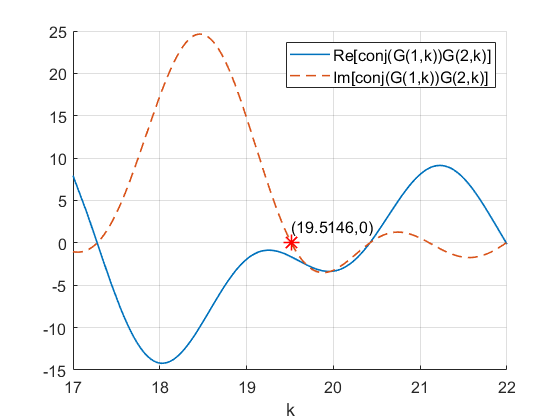}\caption{Left:
representation of the vectors $G$ normalized rotating counter-clockwise for
several values of $k\in[19.31,19.53]$. Right: representation of the functions
$k\rightarrow\text{Re}\left(  \overline{G\left(  z_{b},k\right)  }G\left(
z_{a},k\right)  \right)  $ and $k\rightarrow\text{Im}\left(  \overline
{G\left(  z_{b},k\right)  }G\left(  z_{a},k\right)  \right)  $ with $z_{a} =
2$ and $z_{b}=1$. }%
\label{fig:G}%
\end{figure}

In Figure \ref{fig:G} we illustrate how the alignment between vectors
$G(z_{a},k_{0})$ and $G(z_{b},k_{0})$ takes place.

\begin{remark}
We notice that Lemma \ref{ZeroesLimitFunct} expresses the fact that for large
values of $z$ the function $G\left(  z,k\right)  $ can be approximated as
$G_{1}\left(  z,k\right)  .$ This function is proportional to the product of
to complex numbers which are obtained rotating at different angular speeds
around a circle with radius one around the point $1$ of the complex plane.
Actually one of the numbers (namely $\zeta_{2}$) rotates much more slowly than
the other. Therefore, after placing $\zeta_{2,a},\ \zeta_{2,b}$ at convenient
positions, the whole problem reduces to placing $\zeta_{1,a},\ \zeta_{1,b}$
also at the correct places in order to obtain that $G\left(  z_{a},k\right)
,\ G\left(  z_{b},k\right)  $ point in opposite directions, something that is
possible due to the fact that the points $\zeta_{1,a},\ \zeta_{1,b}$ rotate
around $1$ at different angular speeds. A key point of the argument is to
choose $\zeta_{1,a},\ \zeta_{2,a},\ \zeta_{1,b}$ and $\zeta_{2,b}$ in such a
way that their modulus is bounded from below. This ensures that $\left\vert
G_{1}\left(  z_{a},k\right)  \right\vert $ and$\ \left\vert G_{1}\left(
z_{b},k\right)  \right\vert $ are bounded from below. This fact, combined with
(\ref{A2E0}) allows to treat $G\left(  z_{a},k\right)  $ and$\ G\left(
z_{b},k\right)  $ as perturbations of $G_{1}\left(  z_{a},k\right)  $ and
$G_{1}\left(  z_{b},k\right)  $ respectively.
\end{remark}

\medskip

In Figure \ref{fig:G} (right) we plot the functions 
$$k\rightarrow
\text{Re}\left(  \overline{G\left(  z_{b},k\right)  }G\left(  z_{a},k\right)
\right) ,\ k\rightarrow\text{Im}\left(  \overline{G\left(  z_{b},k\right)
}G\left(  z_{a},k\right)  \right)  $$
 for the values $z_{b}=1,\ z_{a}=2.$ These
functions allow to identify the alignment of the vectors $G\left(
z_{a},k\right)  $ and $G\left(  z_{b},k\right)  $. Indeed, given two complex
numbers $Z_{1},Z_{2}\in\mathbb{C}$ we can identify them with vectors
$V_{1},V_{2}\in\mathbb{R}^{2}$ given by $V_{j}=\left(  \text{Re}\left(
Z_{j}\right)  ,\text{Im}\left(  Z_{j}\right)  \right)  ,\ j=1,2.$ We have the
following identity
\[
\bar{Z}_{1}Z_{2}=V_{1}\cdot V_{2}+i\det\left(  V_{1},V_{2}\right).
\]

Therefore, the two vectors associated to the complex numbers $Z_{1},\ Z_{2}$
are parallel if $\det\left(  V_{1},V_{2}\right)  =\text{Im}\left(  \bar{Z}%
_{1}Z_{2}\right)  =0.$ Moreover, they point in opposite directions if in
addition $V_{1}\cdot V_{2}=\text{Re}\left(  \bar{Z}_{1}Z_{2}\right)  <0.$

Figure \ref{fig:G} (left) shows the existence in the range $k\in\left[
19,20\right]  $ of one value of $k$ such that $G\left(  z_{b},k\right)
=-cG\left(  z_{a},k\right)  $ for some $c<0,\ c\in\mathbb{R}$.

We now come back to the proof of Theorem \ref{LinearzPb}. \medskip

\begin{proofof}
[Proof of Theorem \ref{LinearzPb}]
$\Psi\left(  k;W_{0}\right)  $ defined by means of (\ref{S3E3a}), can be
written, using (\ref{S3E1}), as
\begin{align*}
\Psi\left(  k;W_{0}\right)   &  =e^{-ikX}\int_{-\infty}^{X}dY\int
_{X+\log\left(  1-e^{Y-X}\right)  }^{\infty}dZ\left[  e^{\frac{1}{2}\left(
Y-Z\right)  }W_{0}\left(  Y-Z\right)  \right]  \left(  e^{ikY}+e^{ikZ}\right)
\\
&  =\int_{-\infty}^{0}dY\int_{\log\left(  1-e^{Y}\right)  }^{\infty}dZ\left[
e^{\frac{1}{2}\left(  Y-Z\right)  }W_{0}\left(  Y-Z\right)  \right]  \left(
e^{ikY}+e^{ikZ}\right)  \ .
\end{align*}

It is readily seen that this function is well defined for $k\in\mathbb{R}$ if
$W_{0}$ satisfies (\ref{S3E2}) and (\ref{S1E7}) holds. The fact that
$\Psi\left(  -k;W_{0}\right)  =\overline{ \Psi\left(  k;W_{0}\right)  }$
follows inmediately from the fact that $W_{0}$ is real in $\mathbb{R}$.

In order to simplify the formula for $\Psi\left(  k;W_{0}\right)  $ we use the
change of variables $Z-Y=z,\ dZ=dz.$ Then, using also that $W_{0}\left(
X\right)  =W_{0}\left(  -X\right)  $ we obtain
\[
\Psi\left(  k;W_{0}\right)  =\int_{-\infty}^{0}e^{ikY}dY\int_{\log\left(
1-e^{Y}\right)  -Y}^{\infty}e^{-\frac{z}{2}}W_{0}\left(  z\right)  \left(
1+e^{ikz}\right)  dz \ .
\]

On the other hand, applying Fubini's Theorem and using that if $z\geq
\log\left(  1-e^{Y}\right)  -Y$ if and only if $Y\geq-\log\left(
e^{z}+1\right)  $, we obtain
\begin{equation}
\Psi\left(  k;W_{0}\right)  =\int_{-\infty}^{\infty}e^{-\frac{z}{2}}%
W_{0}\left(  z\right)  \left(  1+e^{ikz}\right)  dz\int_{-\log\left(
e^{z}+1\right)  }^{0}e^{ikY}dY. \label{S3E3c}%
\end{equation}
Thus, computing the integral in $Y$ we obtain
\begin{equation}
\Psi\left(  k;W_{0}\right)  =\frac{1}{ik}\int_{-\infty}^{\infty}e^{-\frac
{z}{2}}W_{0}\left(  z\right)  \left(  1+e^{ikz}\right)  \left[  1-\frac
{1}{\left(  e^{z}+1\right)  ^{ik}}\right]  dz. \label{S3E4a}%
\end{equation}
Using the symmetry $W_{0}\left(  z\right)  =W_{0}\left(  -z\right)  $ we can
rewrite $\Psi\left(  k;W_{0}\right)  $ as
\begin{equation}
\Psi\left(  k;W_{0}\right)  =\frac{1}{ik}\int_{0}^{\infty}W_{0}\left(
z\right)  G\left(  z,k\right)  dz \label{S3E4}%
\end{equation}
where $G\left(  z,k\right)  $ is as in (\ref{S3E5}). It is readily seen that
for each fixed $k\in\mathbb{R}$ we have $G\left(  z,k\right)  =O\left(
e^{-\frac{\left\vert z\right\vert }{2}}\right)  $ as $\left\vert z\right\vert
\rightarrow\infty,\ z\in\mathbb{R}$. Then, if $W_{0}$ satisfies (\ref{S3E2})
we obtain that the integral in (\ref{S3E4}) is well defined due to (\ref{S1E7}).

Due to (\ref{S3E3a}) the problem has been reduced to finding a function
$W_{0}$, with the properties $(i),(ii),(iii)$ stated in the Theorem, such that
the corresponding function $\Psi\left(  k;W_{0}\right)  $ has a zero
$k\in\mathbb{R}\setminus\left\{  0\right\}  .$

In order to show that such a function $W_{0}$ exists we will use the
continuity properties of the function $\Psi\left(  k;W_{0}\right)  $ with
respect to $W_{0}$ if $W_{0}$ experiences small changes in the weak topology
of measures. We will take $W_{0}$ as an even perturbation of the following
combination of Dirac measures
\begin{equation}
W_{0,s}\left(  z\right)  =a\delta\left(  z-z_{a}\right)  +b\delta\left(
z-z_{b}\right)  \ \ ,\ z>0, \ \label{A1E9}%
\end{equation}
where $a>0, b>0, z_{a}>0,\ z_{b}>0.$ Using (\ref{S3E4}) we obtain%
\begin{equation}
\Psi\left(  k;W_{0,s}\right)  =\frac{1}{ik}\left[  aG\left(  z_{a},k\right)
+bG\left(  z_{b},k\right)  \right]. \label{A2E7}%
\end{equation}


In order to obtain a function $W_{0}$ with the regularity and the asymptotic
behaviour stated in the Theorem, we introduce some auxiliary functions
\begin{equation}
\zeta\left(  z\right)  =e^{-z^{2}}\ \ ,\ \ \ \zeta_{\varepsilon}\left(
z\right)  =\frac{1}{\varepsilon}\zeta\left(  \frac{z}{\varepsilon}\right)
\ \ ,\ \ \varepsilon>0,\ \ \ z \in{\mathbb{R}}. \label{S3E6}%
\end{equation}

We define
\begin{equation}
W_{0,1}^{\varepsilon}(z)=a\zeta_{\varepsilon}\left(  z-z_{a}\right)
+b\zeta_{\varepsilon}\left(  z-z_{b}\right)  ,\ z>0 \label{S3E6a}%
\end{equation}
where $a,\ b,\ z_{a},\ z_{b}$ are as in (\ref{A1E9}). Notice that
$W_{0,1}^{\varepsilon}\rightharpoonup W_{0,s}$ as $\varepsilon\rightarrow0$ in
the weak topology of $\mathcal{M}_{+}\left(  \mathbb{R}\right)  .$ We define
$\Psi\left(  k;W_{0,1}^{\varepsilon}\right)  $ using (\ref{S3E4}). It readily
follows that $\lim_{\varepsilon\rightarrow0}\Psi\left(  k;W_{0,1}%
^{\varepsilon}\right)  =\Psi\left(  k;W_{0,s}\right)  $ uniformly in compact
sets of $k.$

On the other hand we define
\begin{equation}
W_{0,2}^{\varepsilon}\left(  z\right)  =\left[  1-\zeta\left(  \varepsilon
z\right)  \right]  \exp\left(  \left(  \frac{\gamma}{2}+p\right)  \sqrt
{z^{2}+1}\right)  ,\ z>0. \label{S3E7}%
\end{equation}
Note that $W_{0,1}^{\varepsilon}$ and $W_{0,2}^{\varepsilon}$ are extended to
the whole real line by means of
\[
\ W_{0,1}^{\varepsilon}\left(  z\right)  = W_{0,1}^{\varepsilon}\left(
-z\right)  , \ \ W_{0,2}^{\varepsilon}\left(  z\right)  = W_{0,2}%
^{\varepsilon}\left(  -z\right)  .
\]

Notice that we can define $\Psi\left(  k;W_{0,2}^{\varepsilon}\right)  $ using
(\ref{S3E4}) since the integral there is convergent for the function
$W_{0,2}^{\varepsilon}$ for each $k\in\mathbb{R}$. Notice that since
$\left\vert G\left(  z,k\right)  \right\vert \leq C_{L}e^{-\frac{\left\vert
z\right\vert }{2}}$ for $\left\vert k\right\vert \leq L$ we have that
$\Psi\left(  k;W_{0,2}^{\varepsilon}\right)  $ converges uniformly to zero as
$\varepsilon\rightarrow0$ uniformly in compact sets of $k$ due to
(\ref{S1E7}). We now define
\begin{equation}
W_{0}^{\varepsilon}=W_{0,1}^{\varepsilon}+W_{0,2}^{\varepsilon} \ .
\label{S3E7a}%
\end{equation}
Then $\Psi\left(  k;W_{0}^{\varepsilon}\right)  =\Psi\left(  k;W_{0,1}%
^{\varepsilon}\right)  +\Psi\left(  k;W_{0,2}^{\varepsilon}\right)  .$ We then
have
\begin{equation}
\lim_{\varepsilon\rightarrow0}\Psi\left(  k;W_{0}^{\varepsilon}\right)
=\Psi\left(  k;W_{0,s}\right)  \label{S3E7b}%
\end{equation}
uniformly in compact sets of $k.$ Moreover, we can write
\begin{equation}
\Psi\left(  k;W_{0,1}^{\varepsilon}\right)  =a\Psi\left(  k;\zeta
_{\varepsilon}\left(  \cdot-z_{a}\right)  \right)  +b\Psi\left(
k;\zeta_{\varepsilon}\left(  \cdot-z_{b}\right)  \right)  \label{A2E8a}%
\end{equation}
and we have also
\begin{equation}
\lim_{\varepsilon\rightarrow0}\Psi\left(  k;\zeta_{\varepsilon}\left(
\cdot-z_{a}\right)  \right)  =G\left(  z_{a},k\right)  \ \ ,\ \ \ \lim
_{\varepsilon\rightarrow0}\Psi\left(  k;\zeta_{\varepsilon}\left(  \cdot
-z_{b}\right)  \right)  =G\left(  z_{b},k\right)  \label{A2E8}%
\end{equation}
uniformly in compact sets of $k.$

We claim that choosing $\varepsilon>0$ we can find $k_{\ast}\in\mathbb{R}$ as
well as $a>0$ and $b>0$ such that $\Psi\left(  k_{\ast};W_{0}^{\varepsilon
}\right)  =0.$ Notice that the function $W_{0}^{\varepsilon}$ depends on $a,b$
although we do not write this dependence explicitly. We write $\Psi\left(
k;\zeta_{\varepsilon}\left(  \cdot-z_{a}\right)  \right)  ,\ \Psi\left(
k;\zeta_{\varepsilon}\left(  \cdot-z_{b}\right)  \right)  $ in polar
coordinates%
\begin{equation}
\Psi\left(  k;\zeta_{\varepsilon}\left(  \cdot-z_{a}\right)  \right)
=R_{a}e^{i\theta_{a}}\ \ ,\ \ \ \Psi\left(  k;\zeta_{\varepsilon}\left(
\cdot-z_{b}\right)  \right)  =R_{b}e^{i\theta_{b}} \label{A3E1}%
\end{equation}
where $R_{a}>0,\ R_{b}>0$ as well as $\theta_{a},\ \theta_{b}$ are functions
of $k.$ Moreover (\ref{A2E8}) combined with Lemma \ref{ZeroesLimitFunct} (in
particular (\ref{S3E5b})-(\ref{S3E5e})) imply that for $\varepsilon>0$
sufficiently small there exist at least one value of $k_{0}\in\mathbb{R}$,
$k_{0}>0$ and some $\delta_{1}>0,\ \sigma>0$ independent of $\varepsilon$ such
that%
\begin{align}
\pi-2\sigma &  \leq\theta_{b}-\theta_{a}\leq\pi-2\sigma\ \ \text{if }%
k=k_{0}-\delta_{1},\label{A2E9a}\\
\pi+\sigma &  \leq\theta_{b}-\theta_{a}\leq\pi+2\sigma\ \ \text{if }%
k=k_{0}+\delta_{1}, \label{A2E9b}%
\end{align}%
\begin{equation}
\frac{1}{L}\leq R_{a}\leq L\ ,\ \frac{1}{L}\leq R_{b}\leq L\ \ \text{if\ }%
k\in\left[  k_{0}-\delta_{1},k_{0}+\delta_{1}\right]  \label{A2E9c}%
\end{equation}
for some $L>0$ independent of $\varepsilon.$

We now select $a,b$ in (\ref{A2E8a}) as follows
\[
a=\frac{1}{R_{a}}\ \ ,\ \ b=\frac{\sigma}{R_{b}},
\]
where $\sigma>0$ is a numerical constant to be determined. Notice that
$R_{a},\ R_{b}$ are functions of $k$ and therefore $a,b$ are also functions.
Due to (\ref{A2E9c}) we have that $\frac{1}{L}\leq a\leq L,\ \frac{\sigma}%
{L}\leq a\leq\sigma L.$ Using (\ref{A2E8a}), (\ref{A3E1}) and our choice of
$a,b $ we obtain%
\[
\Psi\left(  k;W_{0}^{\varepsilon}\right)  =e^{i\theta_{a}}+\sigma
e^{i\theta_{b}}+\Psi\left(  k;W_{0,2}^{\varepsilon}\right)  .
\]

We assume that $\left\vert \sigma-1\right\vert \leq\frac{1}{2}.$ Using the fact
that $\Psi\left(  k;W_{0,2}^{\varepsilon}\right)  $ tends to zero as
$\varepsilon\rightarrow0$ we obtain, combining \ continuity argument with
(\ref{A2E9a}), (\ref{A2E9b}) that if $\varepsilon$ is sufficiently small there
exists $k_{\ast}\in\left(  k_{0}-\delta_{1},k_{0}+\delta_{1}\right)  $ such
that $\left(  e^{i\theta_{a}}+\Psi\left(  k_{\ast};W_{0,2}^{\varepsilon
}\right)  \right)  =-\lambda_0 e^{i\theta_{b}}$ for some $\lambda_0>0$. Notice
that $\lambda_0$ is close to $1$ if $\varepsilon$ is sufficiently small.
Choosing then $\sigma=\lambda_0$ we obtain that $\Psi\left(  k_{\ast}%
;W_{0}^{\varepsilon}\right)  =0.$ Using also $\Psi\left(  -k;W_{0}\right)
=\overline{\Psi\left(  k;W_{0}\right)  }$ we obtain that (\ref{S3E3}) holds
with $W_{0}=W_{0}^{\varepsilon}.$ Notice that $W_{0}$ is analytic in the
domain $D_{\beta}.$

It only remains to prove the asymptotic formula (\ref{S3E3b}). Since this will
imply that $\Psi\left(  k;W_{0}\right)  $ will be different from zero for
large $\left\vert k\right\vert .$ In particular, this implies that we can
define $k_{\ast}$ as the positive root of $\Psi\left(  k;W_{0}\right)  $ in
the real line with the largest value of $\left\vert k\right\vert .$ Using
(\ref{S3E4a}) we obtain%
\[
\Psi\left(  k;W_{0}\right)  =\frac{1}{ik}\left[  G_{1}\left(  k\right)
+G_{2}\left(  k\right)  \right]
\]
where%
\begin{align*}
G_{1}\left(  k\right)  =\int_{-\infty}^{\infty}e^{-\frac{z}{2}}W_{0}\left(
z\right)  \left[  1-\frac{1}{\left(  e^{z}+1\right)  ^{ik}}\right]  dz ,
\ G_{2}\left(  k\right)  =\int_{-\infty}^{\infty}e^{-\frac{z}{2}}%
W_{0}\left(  z\right)  e^{ikz}\left[  1-\frac{1}{\left(  e^{z}+1\right)
^{ik}}\right]  dz .
\end{align*}
We now write
\begin{align*}
&  G_{1}\left(  k\right)  =G_{1,1}\left(  k\right)  +G_{1,2}\left(  k\right)
\\
&  G_{1,1}\left(  k\right)  =\int_{-\infty}^{0}e^{-\frac{z}{2}}W_{0}\left(
z\right)  \left[  1-\frac{1}{\left(  e^{z}+1\right)  ^{ik}}\right]  dz,
\ G_{1,2}\left(  k\right)  =\int_{0}^{\infty}e^{-\frac{z}{2}}W_{0}\left(
z\right)  \left[  1-\frac{1}{\left(  e^{z}+1\right)  ^{ik}}\right]  dz
\end{align*}
and
\begin{align*}
&  G_{2}\left(  k\right)  =G_{2,1}\left(  k\right)  +G_{2,2}\left(  k\right)
\\
&  G_{2,1}\left(  k\right)  =\int_{-\infty}^{0}e^{-\frac{z}{2}}W_{0}\left(
z\right)  e^{ikz}\left[  1-\frac{1}{\left(  e^{z}+1\right)  ^{ik}}\right]  dz, \\
&  \ \ \ G_{2,2}\left(  k\right)  =\int_{0}^{\infty}e^{-\frac{z}{2}}%
W_{0}\left(  z\right)  e^{ikz}\left[  1-\frac{1}{\left(  e^{z}+1\right)
^{ik}}\right]  dz .
\end{align*}
Using the change of variables $\xi=\log\left(  e^{z}+1\right)  $ we readily
obtain, using Riemann-Lebesgue Lemma, as well as (\ref{S3E2}), that
\[
\int_{0}^{\infty}e^{-\frac{z}{2}}W_{0}\left(  z\right)  \frac{dz}{\left(
e^{z}+1\right)  ^{ik}}\rightarrow0\text{ as }\left\vert k\right\vert
\rightarrow\infty\ .
\]
Therefore
\begin{equation}
G_{1,2}\left(  k\right)  \rightarrow\int_{0}^{\infty}e^{-\frac{z}{2}}%
W_{0}\left(  z\right)  dz\text{ as }\left\vert k\right\vert \rightarrow
\infty\ . \label{S5E6}%
\end{equation}
Similarly, we can write $G_{2,2}\left(  k\right)  $ as
\[
G_{2,2}\left(  k\right)  =\int_{0}^{\infty}e^{-\frac{z}{2}}W_{0}\left(
z\right)  e^{ikz}dz-\int_{0}^{\infty}e^{-\frac{z}{2}}W_{0}\left(  z\right)
\left(  \frac{e^{z}}{e^{z}+1}\right)  ^{ik}dz \ .
\]
The first integral on the right converges to zero as $\left\vert k\right\vert
\rightarrow\infty$ due to Riemann-Lebesgue. On the other hand, using the
change of variables $\xi=\log\left(  \frac{e^{z}}{e^{z}+1}\right)  $ in the
second integral on the right-hand side we obtain also that the resulting
integral converges to zero as $\left\vert k\right\vert \rightarrow\infty,$
using again Riemann-Lebesgue. Thus
\begin{equation}
G_{2,2}\left(  k\right)  \rightarrow0\text{ as }\left\vert k\right\vert
\rightarrow\infty.\label{S5E7}%
\end{equation}

It remains to study the asymptotics of $G_{1,1}\left(  k\right)  $ and
$G_{2,1}\left(  k\right)  $ as $\left\vert k\right\vert \rightarrow\infty.$
Using the change of variables $z\rightarrow-z,$ we can rewrite $G_{1,1}%
,G_{2,1}$ as
\begin{align*}
G_{1,1}\left(  k\right)   &  =\int_{0}^{\infty}e^{\frac{z}{2}}W_{0}\left(
z\right)  \left[  1-\frac{1}{\left(  1+e^{-z}\right)  ^{ik}}\right]  dz\\
G_{2,1}\left(  k\right)   &  =\int_{0}^{\infty}e^{\frac{z}{2}}W_{0}\left(
z\right)  e^{-ikz}\left[  1-\frac{1}{\left(  1+e^{-z}\right)  ^{ik}}\right]
dz \ 
\end{align*}
where we used also that $W_{0}\left(  z\right)  =W_{0}\left(  -z\right)  .$
Setting now
\begin{equation}
\label{def:H}H\left(  z;k\right)  =\int_{0}^{z}e^{\frac{\xi}{2}}W_{0}\left(
\xi\right)  e^{-ik\xi}d\xi
\end{equation}
we can rewrite $G_{1,1}, G_{2,1}$ as
\begin{align}
G_{1,1}\left(  k\right)   &  =\int_{0}^{\infty}\frac{\partial H\left(
z;0\right)  }{\partial z}\left[  1-\frac{1}{\left(  1+e^{-z}\right)  ^{ik}%
}\right]  dz\label{S5E7a}\\
G_{2,1}\left(  k\right)   &  =\int_{0}^{\infty}\frac{\partial H\left(
z;k\right)  }{\partial z}\left[  1-\frac{1}{\left(  1+e^{-z}\right)  ^{ik}%
}\right]  dz . \label{S5E7b}%
\end{align}

Integrating by parts in \eqref{S5E7a}, i.e. the formula of $G_{1,1}\left(
k\right)  $, we obtain
\[
G_{1,1}\left(  k\right)  =ik\int_{0}^{\infty}\left[  H\left(  z;0\right)
\frac{e^{-z}}{1+e^{-z}}\right]  \frac{dz}{\left(  1+e^{-z}\right)  ^{ik}}.
\]

Using (\ref{S3E2}) we obtain that, since $1+\gamma+2p\geq1>0,$ the following
asymptotics holds
\[
H\left(  z;0\right)  \frac{e^{-z}}{1+e^{-z}}\sim\frac{2}{1+\gamma
+2p}e^{\left(  \frac{\gamma}{2}+p-\frac{1}{2}\right)  z}\left[  1+O\left(
e^{-\kappa\left\vert z\right\vert }\right)  \right]  \text{ as }%
\text{Re}\left(  z\right)  \rightarrow\infty,\ z\in D_{\beta}%
\]
with $\kappa>0.$ Using the change of variables $X=\log\left(  1+e^{-z}\right)
$ (and hence $z=\log\left(  \frac{1}{e^{X}-1}\right)  $, $dz=-\frac{e^{X}%
}{e^{X}-1}dX$) we can then write
\begin{equation}
G_{1,1}\left(  k\right)  =ik\int_{0}^{\log\left(  2\right)  }\frac
{e^{X}F\left(  X\right)  }{e^{X}-1}e^{-ikX}dX \label{eq:G_11}%
\end{equation}
where $F\left(  \log\left(  1+e^{-z}\right)  \right)  =\frac{e^{-z}H\left(
z;0\right)  }{1+e^{-z}}.$ The function $F$ is analytic in a wedge around the
interval $\left[  0,\log\left(  2\right)  \right]  $ and it satisfies
\begin{equation}
F\left(  X\right)  =\frac{2}{1+\gamma+2p}\left(  X\right)  ^{\frac{1}%
{2}-\left(  \frac{\gamma}{2}+p\right)  }\left[  1+O\left(  \left\vert
X\right\vert ^{\kappa}\right)  \right]  \text{ as }\left\vert X\right\vert
\rightarrow0\ \ \text{if\ }1+\gamma+2p>0 \label{A3E2}%
\end{equation}
with $X$ contained in the portion of a cone $\left\{  \left\vert
\text{Im}\left(  X\right)  \right\vert \leq\theta\left\vert \text{Re}\left(
X\right)  \right\vert :\left\vert X\right\vert \leq\delta_{0}\right\}  $ with
$\theta>0$ and $\delta_{0}>0$ small. Similar asymptotic formulas for the
derivatives of $F\left(  X\right)  $ can be obtained differentiating formally
both sides of (\ref{A3E2}). Using this approximation for $F\left(  X\right)  $
in (\ref{eq:G_11}) and using the change of variables $kX\rightarrow Y$ into
the integral, we obtain the following asymptotic behaviour of $G_{1,1}\left(
k\right)  $ as $|k|\rightarrow\infty$:
\begin{align}
G_{1,1}\left(  k\right)   &  =a\,\text{sgn}\left(  k\right)  e^{i\frac
{\pi\text{sgn}\left(  k\right)  }{2}\left(  \frac{\gamma}{2}+p-\frac{1}%
{2}\right)  }\left\vert k\right\vert ^{\frac{1}{2}+\left(  \frac{\gamma}%
{2}+p\right)  }\left[  1+O\left(  \frac{1}{\left\vert k\right\vert ^{\kappa}%
}\right)  \right]  \text{ as }\left\vert k\right\vert \rightarrow
\infty\label{S5E8}\\
\text{with }a  &  =\frac{2i}{1+\gamma+2p}\Gamma\left(  \frac{1}{2}-\left(
\frac{\gamma}{2}+p\right)  \right)  . \label{S5E8bis}%
\end{align}
Indeed, from \eqref{eq:G_11}, using \eqref{A3E2} we have
\begin{align*}
ik\int_{0}^{\log\left(  2\right)  }\frac{e^{X}F\left(  X\right)  }{e^{X}%
-1}e^{-ikX}dX  &  \sim ik\int_{0}^{\infty}\frac{F\left(  X\right)  }%
{X}e^{-ikX}dX\\
&  \sim\frac{2}{1+\gamma+2p}\left(  ik\right)  \int_{0}^{\infty}\left(
X\right)  ^{-\frac{1}{2}-\left(  \frac{\gamma}{2}+p\right)  }e^{-ikX}dX.
\end{align*}
Moreover, we can rewrite the integral on the right hand side of the equation
above as
\begin{align*}
\int_{0}^{\infty}\left(  X\right)  ^{-\frac{1}{2}-\left(  \frac{\gamma}%
{2}+p\right)  }e^{-ikX}dX  &  =\dots=\left\vert k\right\vert ^{\left(
\frac{\gamma}{2}+p\right)  -\frac{1}{2}}\left(  -i\text{sgn}\left(  k\right)
\right)  \int_{\mathbb{R}_{+}}\left(  -i\text{sgn}\left(  k\right)  Z\right)
^{-\frac{1}{2}-\left(  \frac{\gamma}{2}+p\right)  }e^{-Z}dZ\\
&  =\left\vert k\right\vert ^{\left(  \frac{\gamma}{2}+p\right)  -\frac{1}{2}%
}e^{i\frac{\pi\text{sgn}\left(  k\right)  }{2}\left(  \frac{\gamma}{2}%
+p-\frac{1}{2}\right)  }\Gamma\left(  \frac{1}{2}-\left(  \frac{\gamma}%
{2}+p\right)  \right)  .
\end{align*}
Therefore, we end up with
\begin{align*}
ik\int_{0}^{\log\left(  2\right)  }\frac{e^{X}F\left(  X\right)  }{e^{X}%
-1}e^{-ikX}dX  &  \sim\frac{2}{1+\gamma+2p}\left(  ik\right)  \left\vert
k\right\vert ^{\left(  \frac{\gamma}{2}+p\right)  -\frac{1}{2}}e^{i\frac
{\pi\text{sgn}\left(  k\right)  }{2}\left(  \frac{\gamma}{2}+p-\frac{1}%
{2}\right)  }\Gamma\left(  \frac{1}{2}-\left(  \frac{\gamma}{2}+p\right)
\right) \\
&  =a\,\text{sgn}\left(  k\right)  e^{i\frac{\pi\text{sgn}\left(  k\right)
}{2}\left(  \frac{\gamma}{2}+p-\frac{1}{2}\right)  }\left\vert k\right\vert
^{\frac{1}{2}+\left(  \frac{\gamma}{2}+p\right)  }%
\end{align*}
with $a$ as in \eqref{S5E8bis}.

\medskip

It now remains to estimate the contribution of $G_{2,1}\left(  k\right)  .$ To
this end we integrate by parts in (\ref{S5E7b}) to obtain
\[
G_{2,1}\left(  k\right)  =\int_{0}^{\infty}\frac{\partial H\left(  z;k\right)
}{\partial z}\left[  1-\frac{1}{\left(  1+e^{-z}\right)  ^{ik}}\right]
dz=ik\int_{0}^{\infty}H\left(  z;k\right)  \frac{e^{-z}}{1+e^{-z}}\frac
{1}{\left(  1+e^{-z}\right)  ^{ik}}dz\ .
\]

We can estimate $H\left(  z;k\right)  $ using again integration by parts.
Then
\begin{align}
H\left(  z;k\right)  & =\int_{\left[  0,z\right]  }e^{\frac{\xi}{2}}W_{0}\left(
\xi\right)  e^{-ik\xi}d\xi \\
&=\frac{1}{ik}\left[  W_{0}\left(  0\right)
-e^{\frac{z}{2}}W_{0}\left(  z\right)  e^{-ikz}\right]  +\frac{1}{ik}%
\int_{\left[  0,z\right]  }\frac{d}{d\xi}\left(  e^{\frac{\xi}{2}}W_{0}\left(
\xi\right)  \right)  e^{-ik\xi}d\xi
\end{align}
hence
\[
\left\vert H\left(  z;k\right)  \right\vert \leq\frac{C}{\left\vert
k\right\vert }e^{\left(  \frac{\gamma}{2}+p+\frac{1}{2}\right)  z}.
\]
Then, using that $\gamma+2p<1$ we obtain
\begin{equation}
\left\vert G_{2,1}\left(  k\right)  \right\vert \leq C\text{ for }\left\vert
k\right\vert \geq1. \label{S5E9}%
\end{equation}
Combining now (\ref{S5E6}), (\ref{S5E7}), (\ref{S5E8}), (\ref{S5E9}) we
obtain
\[
\Psi\left(  k;W_{0}\right)  =\frac{1}{ik}\left[  G_{1}\left(  k\right)
+G_{2}\left(  k\right)  \right]  \sim a\text{sgn}\left(  k\right)
e^{i\frac{\pi\text{sgn}\left(  k\right)  }{2}\left(  \frac{\gamma}{2}%
+p-\frac{1}{2}\right)  }\left\vert k\right\vert ^{\frac{1}{2}+\left(
\frac{\gamma}{2}+p\right)  }\text{ as }\left\vert k\right\vert \rightarrow
\infty
\]
whence (\ref{S3E3b}) follows.
\end{proofof}

\bigskip

\textbf{Declaration of interest}
The authors declare that they have no conflict of interest.

\textbf{Acknowledgements}
The authors gratefully acknowledge the support of the Hausdorff Research
Institute for Mathematics (Bonn), through the {\em Junior Trimester
Program on Kinetic Theory\/}, of the CRC 1060 {\em The mathematics of
emergent effects\/} at the University of Bonn funded through the German
Science Foundation (DFG), as well as of the {\em Atmospheric
Mathematics\/} (AtMath) collaboration of the Faculty of Science of
University of Helsinki.
The research has been supported by the Academy of Finland, via
an Academy project (project No. 339228) and the Finnish {\em centre of
excellence in Randomness and STructures\/} (project No. 346306).  The
research of MF has also been partially funded from the ERC Advanced
Grant 741487.
JL and AN would also like to thank the Isaac Newton Institute for
Mathematical
Sciences, Cambridge, for support and hospitality during the programme
{\em Frontiers in kinetic theory: connecting microscopic to macroscopic\/}
(KineCon 2022) where partial work on this paper was undertaken.  This
work was supported by EPSRC grant no EP/R014604/1 and by a grant from
the Simons Foundation.

 \def\adresse{
\begin{description}

\item[M.~A. Ferreira] 
{Department of Mathematics and Statistics, University of Helsinki,\\ P.O. Box 68, FI-00014 Helsingin yliopisto, Finland \\
E-mail:  \texttt{marina.ferreira@helsinki.fi}}\\
ORCID 0000-0001-5446-4845

\item[J. Lukkarinen]{ Department of Mathematics and Statistics, University of Helsinki, \\ P.O. Box 68, FI-00014 Helsingin yliopisto, Finland \\
E-mail: \texttt{jani.lukkarinen@helsinki.fi}}\\
ORCID 0000-0002-8757-1134

\item[A. Nota:] {Department of Information Engineering, Computer Science and Mathematics,\\ University of L'Aquila, 67100 L'Aquila, Italy \\
E-mail: \texttt{alessia.nota@univaq.it}}\\
ORCID 0000-0002-1259-4761

\item[J.~J.~L. Vel\'azquez] { Institute for Applied Mathematics, University of Bonn, \\ Endenicher Allee 60, D-53115 Bonn, Germany\\
E-mail: \texttt{velazquez@iam.uni-bonn.de}}

\end{description}
}

\adresse


\begin{thebibliography}{99}                                                                                               %


\bibitem {A99}D.J. Aldous, Deterministic and Stochastic Models for Coalescence
(Aggregation, Coagulation): A Review of the Mean-Field Theory for
Probabilists, \textit{Bernoulli} \textbf{5}, 3--48 (1999)

\bibitem {BZ}A.M. Balk, V.E. Zakharov, Stability of weak-turbulence
Kolmogorov spectra, \emph{American Mathematical Society Translations},
\textbf{182}(2), 31--82 (1998)

\bibitem {BNV2}M. Bonacini, B. Niethammer, J.J.L. Vel\'{a}zquez, Solutions
with peaks for a coagulation-fragmentation equation. Part I: stability of the
tails. \emph{Communications in Partial Differential Equations}, \textbf{45}
(5), 351--391 (2020)

\bibitem {BNV1}M. Bonacini, B. Niethammer, J.J.L. Vel\'{a}zquez, Solutions
with peaks for a coagulation-fragmentation equation. Part II: Aggregation in
peaks. \emph{Annales de l'Institut Henri Poincar\'e C, Analyse non lin\'eaire}
\textbf{38}(3), 601--646 (2021)

\bibitem {CRZ}C. Connaughton, R. Rajesh, O. Zaboronski, Stationary Kolmogorov
solutions of the Smoluchowski aggregation equation with a source term,
\emph{Phys. Rev. E} 69.6 (2004) 061114

\bibitem {Dub}P.B. Dubovski, Mathematical theory of coagulation, Lecture notes
series, Vol. 23, Seoul National University, Seoul, (1994)

\bibitem {Du}R. Durrett, Random graph dynamics, Vol. 200. No. 7. Cambridge:
Cambridge university press (2007)

\bibitem {FLNV1}M.A. Ferreira, J. Lukkarinen, A. Nota, J.J.L. Vel\'{a}zquez,
Stationary non-equilibrium solutions for coagulation systems. \emph{Arch.
Rational Mech. Anal.} \textbf{240}, 809--875 (2021)

\bibitem {FLNV2}M.A. Ferreira, J. Lukkarinen, A. Nota, J.J.L. Vel\'{a}zquez,
Localization in stationary non-equilibrium solutions for multicomponent
coagulation systems. \emph{Commun. Math. Phys.} \textbf{388} (1), 479--506 (2021)

\bibitem {FLNV3}M.A. Ferreira, J. Lukkarinen, A. Nota, J.J.L. Vel\'{a}zquez,
Multicomponent coagulation systems: existence and non-existence of stationary
non-equilibrium solutions. arXiv: 2103.12763 (2021)

\bibitem {Fried}S.K. Friedlander, Smoke, Dust, and Haze, Oxford University
Press (2000).

\bibitem {GKO}S. Grosskinsky, C. Klingenberg, K. Oelschl\"ager, A rigorous
derivation of Smoluchowski's equation in the moderate limit, \emph{Stoch.
Anal. Appl.}, \textbf{22}, 113--141 (2004)

\bibitem {HR}A. Hammond, F. Rezakhanlou, The kinetic limit of a system of
coagulating Brownian particles, \emph{Arch. Rat. Mech. Anal.} \textbf{185}(1),
1--67 (2007)

\bibitem {H}H. Hayakawa, Irreversible kinetic coagulations in the presence of
a source, \emph{Journal of Physics A: Mathematical and General}, \textbf{20}
(12) L801--L805 (1987)

\bibitem {HNV}M. Herrmann, B. Niethammer, J.J.L. Vel\'{a}zquez, Instabilites
and oscillations in coagulation equations with kernels of homogeneity one,
\emph{Quarterly Appl. Math.}, \textbf{75}(1), 105--130, (2017)

\bibitem {KC}P.L. Krapivsky and C. Connaughton, Driven Brownian coagulation
of polymers, \emph{J. Chem. Phys.} \textbf{136}, 204901 (2012)

\bibitem {LN}R. Lang, X.-X. Nguyen, Smoluchowski's theory of coagulation in
colloids holds rigorously in the Boltzmann-Grad-limit, \emph{Z. Wahrsch. Verw.
Gebiete} \textbf{54}(3), 227--280 (1980)

\bibitem {LNV}P. Lauren\c{c}ot, B. Niethammer, J.J.L. Vel\'{a}zquez,
Oscillatory dynamics in Smoluchowski's coagulation equation with diagonal
kernel, \emph{Kinet. Relat. Models}, \textbf{11}, 933--952 (2018)

\bibitem {McNV}J. B. McLeod, B. Niethammer, and J.J.L. Vel\'{a}zquez,
Asymptotics of self-similar solutions to coagulation equations with product
kernel, \emph{J. Stat. Phys.} \textbf{144}(1), 76--100 (2011)

\bibitem {MSST}S.A. Matveev, A.A. Sorokin, A.P. Smirnov and E.E.
Tyrtyshnikov, Oscillating stationary distributions of nanoclusters in an open
system, \emph{Mathematical and Computer Modelling of Dynamical Systems},
\textbf{26}(6), 562--575 (2020)

\bibitem {MKSTB} S.A. Matveev, P.L. Krapivsky, A.P. Smirnov, 
E.E. Tyrtyshnikov, N.V. Brilliantov,  Oscillations in aggregation-shattering
processes. \emph{Physical Review Letters}, \textbf{119}(26), 260601 (2017)

\bibitem {NPSV21}B. Niethammer, R.L. Pego, A. Schlichting, J.J.L. Vel\'azquez, Oscillations in a Becker-D\" oring model with injection and depletion. To appear in \emph{SIAM J. Appl. Math.} (2022)


\bibitem {NV}B. Niethammer and J.J.L. Vel\'azquez, Oscillatory traveling
wave solutions for coagulation equations \emph{Quart. Appl. Math.},
\textbf{76}(1):153--188, (2018)

\bibitem {NoV}A. Nota, J.J.L. Vel\'azquez, On the Growth of a Particle
Coalescing in a Poisson Distribution of Obstacles, \emph{Commun. Math. Phys.}
\textbf{354}, 957--1013 (2017)

\bibitem {PV20} R.L. Pego,  J.J.L. Vel\'{a}zquez, Temporal oscillations in
Becker-D\"{o}ring equations with atomization. \emph{Nonlinearity},
\textbf{33}(4), 1812 (2020)

\bibitem {Ru87}W. Rudin, Real and Complex Analysis,McGraw-Hill Education, (1987)

\bibitem {TIN}H. Tanaka, S. Inaba, K. Nakazawa, Steady-state size distribution
for the self-similar collision cascade, \emph{Icarus} \textbf{123}(2) (1996)
450-- 455.

\bibitem {YRH}M.R. Yaghouti, F. Rezakhanlou, A. Hammond, Coagulation,
diffusion and the continuous Smoluchowski equation, \emph{Stoch. Proc. Appl.}
\textbf{119}(9), 3042--3080 (2009)

\bibitem {V} M.S. Veshchunov, A new approach to the Brownian coagulation
theory. Journal of Aerosol Science 41.10 (2010): 895-910

\bibitem {Zakh}V.E. Zakharov, V.S. L'vov, G. Falkovich, Kolmogorov spectra of
turbulence I. Wave turbulence. Springer-Verlag Berlin Heidelberg (1992)

\bibitem {ZF}V.E. Zakharov, N.N. Filonenko, Energy spectrum for stochastic
oscillations of a fluid surface, \emph{Dokl. Acad. Nauka SSSR}, \textbf{170}
(1966) 1292--1295, \emph{Sov. Phys. Dokl.} \textbf{11} (1967) 881- 884.
\end{thebibliography}
\end{document}